\documentclass[11pt,twoside,a4paper]{article}

\newcommand{\comment}[1]{}

\usepackage{t1enc} 
\usepackage[dvips]{epsfig}
\usepackage{amssymb}
\usepackage{amsmath}
\usepackage{longtable}
\usepackage{makeidx}
\usepackage{euscript}
\usepackage{enumerate}
\usepackage{fancyhdr}
\usepackage[usenames, dvipsnames]{color}

\DeclareMathAlphabet\scr{U}{scr}{m}{n}
\SetMathAlphabet\scr{bold}{U}{scr}{b}{n}
  \DeclareFontFamily{U}{scr}{\skewchar\font'177}%
  \DeclareFontShape{U}{scr}{m}{n}{<-6>rsfs5<6-8>rsfs7<8->rsfs10}{}%
  \DeclareFontShape{U}{scr}{b}{n}{<-6>rsfs5<6-8>rsfs7<8->rsfs10}{}%

\newcommand{\rr}{\mathbb R}  
\newcommand{\ev}{\mathbb E}

\newcommand{\nn}{\mathbb N}

\newcommand{\pp}{\mathbb P}

\newcommand{\iprod}[1]{\left\langle #1 \right\rangle}

\newtheorem{satz}{Theorem}[section]
\newtheorem{prop}[satz]{Proposition}
\newtheorem{theorem}[satz]{Theorem}

\newtheorem{lemma}[satz]{Lemma}
\newtheorem{cor}[satz]{Corollary}
\newtheorem{assumption}[satz]{Assumption}
\newtheorem{@definition}[satz]{Definition}

\newtheorem{@aufgabe}{Aufgabe}

\newtheorem{@bsp}[satz]{Beispiel}

\newenvironment{bem}{\noindent {\bf Remark.\;}}{}

\newenvironment{mat}{\left (\begin{matrix} } {\end{matrix}\right ) }

\newcommand{\bmat}{\begin{mat}}
\newcommand{\emat}{\end{mat}}
\newcommand{\be}{\begin{enumerate}}
\newcommand{\ee}{\end{enumerate}}
\newcommand{\beq}{\begin{equation}}
\newcommand{\eeq}{\end{equation}}
\newcommand{\bea}{\begin{eqnarray}}
\newcommand{\eea}{\end{eqnarray}}
\newcommand{\beaa}{\begin{eqnarray*}}
\newcommand{\eeaa}{\end{eqnarray*}}

\newcommand{\ep}{\hfill $\Box$}

\newcommand{\eq}{\begin{equation}}
\newcommand{\en}{\end{equation}}
\newcommand{\phibar}{\overline{\Phi}}
\newcommand{\abs}[1]{\left\lvert #1 \right\rvert}
\newcommand{\met}{\mathcal{X}}
\newcommand{\wass}{\mathcal{W}}
\newcommand{\ent}{H}
\newcommand{\lip}{\mathcal{L}}

\newcommand{\newz}{\widetilde{Z}_{(1)}}

\newcommand{\norm}[1]{\left\lVert #1 \right\rVert}

\renewcommand{\epsilon}{\varepsilon}
\renewcommand{\phi}{\varphi}
\renewcommand{\rho}{\varrho}

\title{{\small\bf CONCENTRATION OF MEASURE FOR BROWNIAN PARTICLE SYSTEMS INTERACTING THROUGH THEIR RANKS}}
\author{{\small BY S. PAL\footnote{This research is partially supported by NSF grant DMS-1007563.}}\\{\small\it University of Washington}\\\quad\\{\small AND}\\\quad\\{\small BY M. SHKOLNIKOV\footnote{Research supported in part by NSF grant DMS-0806211.}}\\{\small\it Stanford University}}
\date{}

\begin{document}

\maketitle

\begin{abstract}
We consider a finite or countable collection of one-dimensional Brownian particles whose dynamics at any point in time is determined by their rank in the entire particle system. Using Transportation Cost Inequalities for stochastic processes we provide uniform fluctuation bounds for the ordered particles, their local time of collisions, and various associated statistics over intervals of time. For example, such processes, when exponentiated and rescaled, exhibit power law decay under stationarity; we derive concentration bounds for the empirical estimates of the index of the power law over large intervals of time. A key ingredient in our proofs is a novel upper bound on the Lipschitz constant of the Skorokhod map that transforms a multidimensional Brownian path to a path which is constrained not to leave the positive orthant.      
\end{abstract}

\section{Introduction}

Define the set $I$ as $\{1,\dots,K\}$ for some $K\in\nn$ or as the set of natural numbers $\nn$ and let $\delta_i$, $i\in I$ be a finite or countable collection of real constants. Consider the following system of stochastic differential equations: 
\begin{equation}\label{ranksde}
dX_i(t) = \sum_{j\in I} \delta_j\cdot1_{\left\{X_i(t)=X_{(j)}(t)\right\}} dt + dW_i(t),\quad i\in I.
\end{equation}
Here, $X_{(1)}(t) \le X_{(2)}(t) \le \ldots $ are the coordinates of the process $X_i(t)$, $i\in I$ in the increasing order and $W=(W_i:\;i\in I)$ is a system of jointly independent one-dimensional standard Brownian motions. The equations in (\ref{ranksde}) model the movement of the particles by interacting Brownian motions such that at every time point, if we order the positions of the particles, then the $i$-th ranked particle from the bottom gets a drift $\delta_i$ for every $i\in I$. As time evolves, the Brownian motions switch ranks and drifts, and, hence, their motion is determined by these time-dependent interactions. When $I$ is finite, the existence and uniqueness in law of such processes is a consequence of Girsanov's Theorem (see e.g. Lemma 6 in \cite{pp}). The countable case is subtle and requires constraints on the initial positions of the particles. In particular, for (\ref{ranksde}) to make sense, the number of particles on every interval of the form $(-\infty,x]$ has to be finite at any point in time with probability one. We discuss this issue in more detail later in the text.

Different versions of the particle system in (\ref{ranksde}) have been considered in several recent articles. Among the more recent ones, see Banner, Fernholz and Karatzas \cite{atlasmodel}, Banner and Ghomrasni \cite{BG}, McKean and Shepp \cite{sheppmckean}, Pal and Pitman \cite{pp}, Jourdain and Malrieu \cite{joumal}, Chatterjee and Pal \cite{chatpal,chatpal2}, Ichiba and Karatzas \cite{ik}, Ichiba et al.~\cite{IPBKF}, and Shkolnikov \cite{sh,sh2}. We refer the reader to the above articles for the full list of applications of such processes. 
Related discrete time processes are studied in the context of the Sherrington-Kirkpatrick model of spin glasses by Arguin and Aizenman \cite{AA}, Ruzmaikina and Aizenman \cite{ruzaizenman} and Shkolnikov \cite{shkol}.

Classically, the case when $\delta_i=0$ for all $i\in I$ has been treated by Harris \cite{harris65}, Arratia \cite{arratia83}
and Sznitman \cite{sznitman86},\cite[p. 187]{sznitman91}, and more recently by Swanson \cite{Sw}. In this case one imagines countably many Brownian particles moving under elastic collision, i.e. bouncing off one another when their paths meet. This is the process of ordered particles derived from the system in \eqref{ranksde}. Harris \cite[(7.1)]{harris65} gave an explicit formula for the law of
$B_{(0)}(t)$, the location at time $t$ of the particle which is the leftmost in the infinite system with $\delta_i=\delta$ for all $i\in I$. From this he deduced for $2 \delta = 1$ that 
\eq\label{har0}
{B_{(0)}(t) \over  t ^{1/4} } \stackrel{d}{\rightarrow}
\left(2 \over \pi \right)^{1/4} {1 \over \sqrt{2 \delta}} \, B(1) \mbox{ as } t \rightarrow \infty 
\en
where $B(1)$ is standard Gaussian.
As remarked by Arratia \cite[p. 71]{arratia83}, the conclusion for general $\delta >0$ follows from the case $2 \delta = 1$ by Brownian scaling.
See also De Masi, Ferrari \cite{demasiferrari}, Rost, Vares \cite{rostvares} and Arratia \cite{arratia83} where variants (or generalizations) of \eqref{har0} are proved for a tagged particle in the exclusion process on $\mathbb{Z}$ associated with a simple symmetric random walk.
Harris conjectured that the process $B_{(0)}$ is not Markov, and left open the problem of describing the long-run
behavior of the paths of $B_{(0)}$. As far as we know, these problems are still open.  

The case of distinct drift parameters differs from these classical models in several remarkable ways. For example, if $I=\{1,\dots,K\}$ for some $K\in\nn$ and 
$$\frac{1}{j}\sum_{k=1}^j \delta_k>\frac{1}{K-j}\sum_{k=j+1}^K \delta_k$$
 for all $j=1,\dots,K-1$, then there exists an invariant distribution for the process of gaps between consecutive particles $(X_{(2)}(t)-X_{(1)}(t),\dots,X_{(K)}(t)-X_{(K-1)}(t))$, $t\geq0$ (see Theorem 8 in \cite{pp}) which is not the case if $\delta_1=\dots=\delta_K$. When $I=\nn$, one can consider the so-called \textit{Atlas} model in which $\delta_1=\delta>0$ and $\delta_i=0$ for every $i\neq1$ in $I$. In this case it is shown in \cite{pp} that, if the initial positions of the particles are chosen according to a standard Poisson process of rate $2\delta$, then the joint distribution of the gaps between consecutive particles  $X_{(2)}(t)-X_{(1)}(t),X_{(3)}(t)-X_{(2)}(t),\dots$ is the same for all $t\geq0$. However, the same statement is not true when $\delta_i=\delta$ for all $i\in I$ (see Theorem 4.2 in \cite{ruzaizenman}). Moreover, if we consider the motion of the left-most particle $X_{(1)}(t)$, $t\geq0$ in the Atlas model, no estimates on its growth and fluctuations are known. Our article is a step in the latter direction. Using techniques from the theory of concentration of measures we give estimates on fluctuations of the paths of the distances between ordered particles and associated statistics. 

One such statistic is given by the \textit{market weights}. The latter can be defined for $I=\{1,\dots,K\}$ with an arbitrary $K\in\nn$ and any choice of $\delta_1,\ldots,\delta_K$ by setting 
\eq\label{whatismu}
\mu_i(t)=\frac{e^{ X_{(K-i+1)}(t)}}{\sum_{j=1}^K e^{ X_{(j)}(t)}}, \qquad i=1,\ldots,K
\en
for all $t\geq0$. It is clear that for any fixed $t$ the sequence $\mu_1(t), \ldots, \mu_K(t)$ is a non-increasing sequence of positive numbers that add up to one. These numbers, in econometric models, go by the name of {market weights} and have an interesting history. 

Fernholz in his 2002 book \cite{fern02} introduces solutions of \eqref{ranksde} to model the time dynamics of the logarithmic market capitalizations of different companies in an equity market. In other words, he considers a stock market with $K$ companies whose total worths in stocks are given by exponentials of the one-dimensional components of the solution of equation \eqref{ranksde}. A major objective of his work is to explain the following curious empirical fact. Consider the shape one obtains by plotting $\log \mu_i$ versus $\log i$. This log-log plot is referred to as the \emph{capital distribution curve}. See Figure 1 below (reproduced from \cite{fern02}) which shows the capital distribution curves between 1929 and 1999 for all the major US stock markets (NYSE, AMEX and NASDAQ) combined. Empirically, the left part of the plot exhibits nearly linear decay. This corresponds to the market weights, in decreasing order, displaying power law (or Zipf's law) decay. More strikingly, the slope of the decay is nearly constant over eight decades, something truly remarkable in the volatile world of financial markets.

\begin{figure}[t]
\begin{center}
\includegraphics[width=4 in]{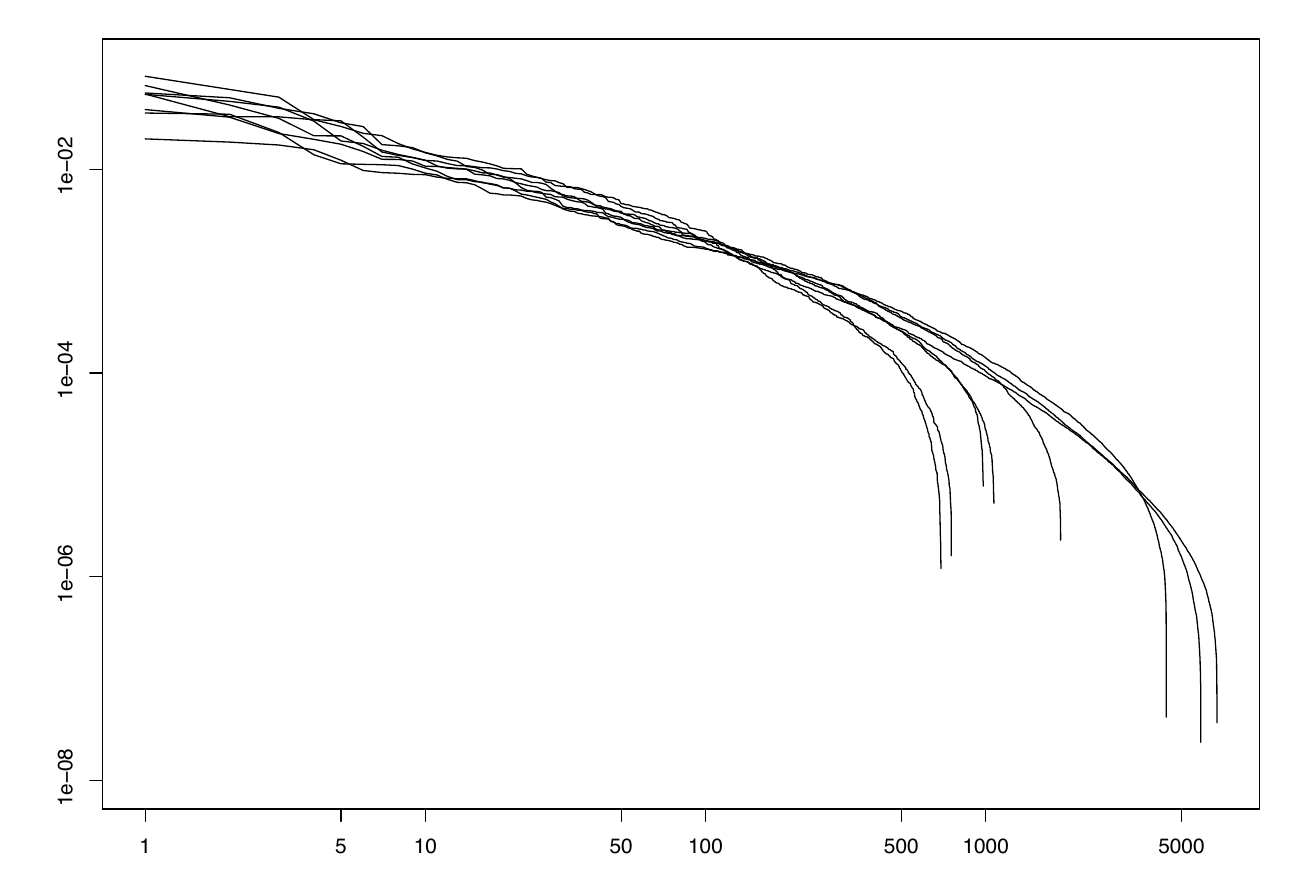}
\caption{Capital distribution curves: 1929-1999}
\end{center}
\end{figure}

In \cite{chatpal} the authors explain the linearity by proving (under suitable assumptions) that the possible limiting stationary laws of the market weights, as $K$ grows to infinity, are given by a subset of the Poisson-Dirichlet family of distributions. The masses of this family have a polynomial decay with size, which corresponds to the linear decay in the log-log plot. However, this does not quite address the temporal stability of the picture above. Our analysis below captures some of its subtleties.  
\bigskip

To analyze the stability of the shape of the capital distribution curve, consider the process of market weights $(\mu_i(t), \; i=1,\dots, K)$, $t\ge 0$ and fix a $J \ll K$. At any point of time $t\geq0$ we introduce the linear regression between the pairs of data $\{ (\log i, \log \mu_{i}(t)),\; i=1,\ldots, J \}$, passing through the first point $(0, \log \mu_1(t))$. In other words, we deal with regressions of the form 
\eq\label{linreg}
\log \mu_1(t) - \log \mu_i(t) = \alpha(t) \log i + \epsilon_i(t), \quad i=2,\ldots, J
\en
for fixed values of $t\geq0$. Clearly, the resulting ordinary least squares estimator for the slope parameter corresponds to the slope of the linear part of the curve at time $t$ in Figure 1. With a minor abuse of notation we will denote this estimator by $\alpha(t)$.

Our objective is to estimate the fluctuations of the above parameter when the spacing process is running close to its stationary law. We choose the Atlas model to have a specific sequence of drifts $\delta_1,\dots,\delta_K$, although other values of the drift parameters can be easily substituted. We consider initial configurations the spacings of which are close to their unique invariant distribution found in \cite{pp} (see Lemma \ref{invatlas} below).

\begin{theorem}\label{slopetheorem}
Let $I=\{1,\dots,K\}$ and consider the Atlas model in which $\delta_1=\delta$ and $\delta_i=0$ for all $i=2,\dots,K$. Define the initial spacings between the particles by
\eq\label{initial_space}
X_{(K)}(0) - X_{(K-i+1)}(0) = \frac{K}{\delta}\log i,\quad i=2,\dots,K.
\en
Assume that $K$ is sufficiently large. Fix a $J < K/15$ and consider the process $\alpha(t)$, $t\geq0$ of ordinary least squares estimators resulting from regressions in \eqref{linreg} for different values of $t$. Then the process $\alpha(t)$, $0\leq t\leq\delta^{-2}K$ satisfies the following concentration of measure property. Let $\overline\alpha= \sup_{0\le s\le \delta^{-2}K} \left[ \alpha(s)   \right]$.

Then there is a constant $m_\alpha$ such that 
\[
\pp\left( \overline\alpha < m_\alpha  \right) \le 1/2 + 17 \sqrt{K} e^{-K/500},
\]
and for all positive $r$ sufficiently large, one has
\[
\pp\left(  \overline\alpha > m_{\alpha} + r \sqrt{K}  \right) \le 2 \exp\left(-\frac{r^2\delta^2}{\mu C_\alpha}  \right).
\]
Hereby, $\mu$ is an absolute positive constant and $C_{\alpha}$ is a positive constant depending on $J$ and given by
\[
C_{\alpha}(J) = J^3\cdot\left(\frac{\sum_{i=1}^{J-1} \log(J!/ i!)}{\sum_{i=2}^J (\log i)^2}\right)^2.
\]
\end{theorem}

\begin{bem} The above result, although novel, does not capture fully the extent of concentration that is seen in the real world data (see Figure 1). We suspect that the reason for this is the empirically observed unequal diffusion coefficients of the ordered particles. See, for example, the discussion in \cite{atlasmodel} which mentions that the diffusion coefficient for the ranked particles decays linearly with increasing rank. Thus, particles at the top are more stable than the average which should lead to a higher concentration. 
\end{bem}
\medskip

\begin{bem}
Another obvious way to improve the bound is to use true stationarity. One can divide up a long time interval in small subintervals of appropriate size, use the above theorem on the small intervals, and take a union bound. However, in this case it is not clear if the Gaussian concentration is preserved when the invariant distribution has Exponential tails. However, we expect that variations of our method and the argument in this remark can lead to Exponential concentration over much larger intervals.  
\end{bem}

\bigskip

Next, we set $I=\nn$ in (\ref{ranksde}) and consider a sequence of drift coefficients $\delta_1,\delta_2,\dots$ which satisfies 
\begin{eqnarray}
\delta_M=\delta_{M+1}=\dots 
\end{eqnarray}
for some $M\in\nn$. Under the assumption that the sequence of initial positions of the particles $X_1(0),X_2(0),\dots$ is deterministic, non-decreasing and such that
\begin{eqnarray}
\liminf_{i\rightarrow\infty}\frac{X_i(0)}{i}>0, 
\end{eqnarray}
the system of stochastic differential equations (\ref{ranksde}) has a unique weak solution (see Proposition 3.1 in \cite{sh}). We let $L_{(i,i+1)}(t)$, $t\geq0$ be the local time process at zero of the process $X_{(i+1)}(t)-X_{(i)}(t)$, $t\geq0$ for each $i\in\nn$ and are interested in the concentration properties of the vector-valued processes $(X_{(2)}(t)-X_{(1)}(t),\dots,X_{(n)}(t)-X_{(n-1)}(t))$, $t\in[0,T]$ and $(L_{(1,2)}(t),\dots,L_{(n-1,n)}(t))$, $t\in[0,T]$ for arbitrary values of $n\in\nn$ and $T\geq0$, which we consider to be fixed from now on. For our main result on the latter we introduce the following assumption on the initial particle configuration $(X_1(0),X_2(0),\dots)$:

{\assumption\label{assumption} The sequence of initial positions of the particles is non-decreasing with probability one, and there exists a deterministic constant $c>0$ such that  
\begin{eqnarray}\label{ic}
X_k(0)-X_N(0)\geq c(k-N)
\end{eqnarray}
holds for all $k\geq N$ almost surely, where $N=\max(n,M)$.\\\quad\\}
We note that Assumption \ref{assumption} ensures the existence of a unique weak solution for the system (\ref{ranksde}) by Proposition 3.1 in \cite{sh}. 

\bigskip

To state our second main theorem we define the norm 
\begin{eqnarray}\label{T2norm}
\|f\|_{T,2}=\Big(\frac{1}{n-1}\sum_{i=1}^{n-1} \sup_{0\leq t\leq T} f_i(t)^2\Big)^{\frac{1}{2}} 
\end{eqnarray} 
on the space $C([0,T],\rr^{n-1})$ of continuous $\rr^{n-1}$-valued functions on $[0,T]$ where $f_i$, $i=1,\dots,n-1$ are the component functions of $f$. Our result then reads as follows.
{\theorem\label{thm_infinite} Let the Assumption \ref{assumption} be satisfied with a constant $c>0$. Moreover, let $A$ and $B$ be measurable subsets of $(C([0,T],\rr^{n-1}),\|.\|_{T,2})$ such that
\begin{eqnarray}
&&\pp\Big(((L_{(i,i+1)}(t), i=1,\dots,n-1),\;t\in[0,T])\in A\Big)>\frac{1}{2}, \label{Acond}\\
&&\pp\Big(((X_{(i+1)}(t)-X_{(i)}(t),i=1,\dots,n-1),\;t\in[0,T])\in B\Big)>\frac{1}{2} \label{Bcond}
\end{eqnarray}
and for any $r>0$ set  
\begin{eqnarray}
&&A_r=\{h\in C([0,T],\rr^{n-1})|\;\inf_{\widetilde{h}\in A}\|h-\widetilde{h}\|_{T,2}\leq r\}, \\
&&B_r=\{h\in C([0,T],\rr^{n-1})|\;\inf_{\widetilde{h}\in B}\|h-\widetilde{h}\|_{T,2}\leq r\}. 
\end{eqnarray}
Then there exist positive constants $C_1$, $C_2$, $C_3$, $C_4$ depending on $c$, $\Delta=\max_{j=1,\dots,M-1} |\delta_j-\delta_M|$, $M$, $n$, $T$  and the value of the left-hand side of (\ref{Acond}) and (\ref{Bcond}), respectively, such that for all $r_1\geq C_1$ and $r_2\geq C_2$ it holds
\begin{eqnarray*}
&&\pp\Big(((L_{(i,i+1)}(t),i=1,\dots,n-1),t\in[0,T])\notin A_{r_1}\Big)
\leq C_3\exp\Big(-r_1^{\frac{4}{7}}\cdot\frac{(n-1)^{\frac{2}{7}}c^{\frac{10}{7}}}{14 T}\Big), \label{Arineq}\\
&&\pp\Big(((X_{(i+1)}(t)-X_{(i)}(t),i=1,\dots,n-1),t\in[0,T])\notin B_{r_2}\Big)
\leq C_4\exp\Big(-r_2^{\frac{4}{7}}\cdot\frac{(n-1)^{\frac{2}{7}}c^{\frac{10}{7}}}{25 T}\Big). \label{Brineq}
\end{eqnarray*}}

\begin{bem} The local time of collisions between two consecutive ordered particles, as considered in Theorem \ref{thm_infinite}, is interesting both mathematically and in applications. Mathematically, say in the Atlas model, the local time (compared to the the case $\delta_i=0$, $i\in I$) measures the \textit{push} felt by the various particles due to the drift at the bottom. Its significance in economic models is discussed in Section 13 of \cite{FK}, which also mentions the somewhat surprising fact that these local times can indeed be measured from data. 
\end{bem}

\bigskip

The rest of the article is organized as follows. In the next section we recall some facts about the concentration of measure phenomenon and the Skorokhod problem in the orthant, and provide an upper bound on the Lipschitz constant for the Skorokhod map of interest in Lemma \ref{lemma_lipschitz}. The latter is the key to the proofs of the two main results. Its proof relies on the construction of the Skorokhod map by Harrison and Reiman \cite{hr} and applies to other Skorokhod problems in the orthant as well. In Section \ref{sec_slope} we use a relation between Transportation Cost Inequalities and the concentration of measure phenomenon to complete the proof of Theorem \ref{slopetheorem} and provide the remainder of the proof of Theorem \ref{thm_infinite} in Section \ref{sec_infinite}. 

\section{Preliminaries}

\subsection{Some facts about measure concentration for stochastic processes} 

Suppose $(\met, d)$ is a complete separable metric space equipped with the Borel $\sigma$-algebra. For a Borel subset $A$ of $\met$ and a positive real number $r$, define the $r$-neighborhood of $A$ by
\[
A_r := \left\{ x\in \met:\quad d(x,A)\le r   \right\}.
\]
We say that a probability measure $\mu$ on $(\met, d)$ has the \textit{measure concentration} property if for any Borel subset $A$ with $\mu(A)\geq1/2$ the value of $\mu(A_r)$ is \textit{very close} to one for large values of $r$. This closeness is usually expressed as a sub-gaussian tail in terms of $r$.

The concentration of measure phenomenon has become one of the most important concepts in modern probability theory. For an excellent introduction to this area we refer the reader to the beautiful article by Talagrand \cite{T96}. A considerable effort has been spent by probabilists on identifying distributions that have the measure concentration property. A (somewhat dated) survey can be found in the monograph by Ledoux \cite{L}.

One technique for proving the measure concentration property, originally proposed by Marton \cite{M1}, involves the so-called Transportation Cost Inequalities (TCI) that we describe below (see also Talagrand \cite{T96b}, Dembo \cite{D}, Dembo and Zeitouni \cite{DZ}). Consider, as before, a complete separable metric space $(\met,d)$ endowed with its Borel $\sigma$-algebra. For a real number $p\geq1$ and all probability measures $P$ and $Q$ on the latter space define the $p$-th Wasserstein distance
\[
\wass_p(P, Q)= \inf_{\pi} \left[ Ed\left( X, X' \right)^p \right]^{1/p},
\]
where the infimum is taken over all couplings of a pair of random elements $(X, X')$ such that the marginal law of $X$ is $P$ and that of $X'$ is $Q$. 

Next, fix a particular probability measure $P$. Suppose there is a constant $C > 0$ such that for all probability measures $Q \ll P$ we have
\eq\label{qtci}
\wass_2(P,Q) \le \sqrt{2C\ent\left( Q \mid P \right) }
\en
where $\ent$ refers to the relative entropy $\ent\left( Q \mid P \right)= E^Q \log\left( {dQ}/{dP} \right)$. In this case we say that $P$ satisfies the \textit{Quadratic Transportation Cost Inequality} (QTCI) with the constant $C$. 

A function $f:\met \rightarrow \rr$ will be called Lipschitz if there is a positive constant $\alpha$ for which
\[
\abs{f(x)- f(y)} \le \alpha d(x,y), \quad x,y \in \met.
\]
The smallest such constant $\alpha$ is then referred to as the Lipschitz constant of $f$. We shall call a function $1$-Lipschitz if $\alpha$ can be taken to be one. Let $\lip$ denote the set of all $1$-Lipschitz functions on $(\met,d)$. The (very short) proof of the following theorem can be found in Ledoux \cite[p.~118]{L} and the original article by Marton \cite{M1}.

\begin{theorem}\label{preq}
Suppose that $P$ satisfies the QTCI with constant $C$. Then one has the following concentration estimates for all $r \ge 2\sqrt{2C\log 2}$:
\begin{enumerate}
\item[(i)] For any measurable set $A$ such that $\mu(A) \ge 1/2$, it holds
\eq\label{aconcen}
\mu(A_r) \ge 1 - \exp\left( - r^2/8C \right).
\en
\item[(ii)] For any $f \in \lip$, one has
\eq\label{fconcen}
P\left( x:\; \abs{f(x) - m_f} \ge r \right) \le 2 e^{-r^2/8C},
\en
where $m_f$ is the median of $f$ with respect to $P$.
\end{enumerate}
\end{theorem}

The QTCI is unique in its advantages and is related to the log-Sobolev inequality, hypercontractivity, Poincar\'e inequality, inf-convolution and Hamilton-Jacobi equations. For details, please consult Otto and Vilani \cite{OV}, Bobkov and G\"otze \cite{BG}, and Bobkov, Gentil and Ledoux \cite{BGL}. 

In addition to Theorem \ref{preq}, the following lemma will be useful in the later text. Its (short) proof can be found in Lemma 2.1 of \cite{DGW}.

\begin{lemma}\label{lemmapushforward}
Suppose $\mu$ is a measure on a metric space $(E, d_E)$ that satisfies the QTCI with a constant $C$. Let $(F, d_F)$ be another metric space. If the map $\Psi:(E,d_E) \rightarrow (F, d_F)$ is Lipschitz, i.e.
\[
d_F\left( \Psi(x), \Psi(y) \right) \le \alpha d_E(x,y), \quad x, y \in E
\]
for some constant $\alpha>0$, then $\tilde\mu=\mu\circ \Psi^{-1}$ satisfies the QTCI on $(F, d_F)$ with the constant $C\alpha^2$.
\end{lemma}

In this article we are interested in the choice of $\met=C([0,T], \rr^{K-1})$, the space of continuous maps from the interval $[0,T]$ to $\rr^{K-1}$, where $K$ is the number of interacting particles as before. The latter space is typically referred to as the \textit{path space}. It will be endowed with different variants of the uniform metric which are described below. We shall use Quadratic Transportation Cost Inequalities satisfied by the laws of stochastic processes with continuous paths, especially multidimensional reflected Brownian motions. 

Let us provide a brief review of the literature on measure concentration in path space. Houdr\'e and Privault \cite{HP}, and Nourdin and Viens \cite{NV} use tools from Malliavin Calculus to derive concentration inequalities for functionals on the Wiener space, among other things. Transportation Cost Inequalities on the Wiener space first appeared in Feyel and \"Ust\"unel \cite{FU}. Djellout, Guillin and Wu \cite{DGW} provide a characterization of the $\mathbf{L}^1$-TCI for diffusions. They also prove the QTCI for diffusions with respect to the Cameron-Martin $\mathbf{L}^2$-metric. Several articles in analysis and geometry are also devoted to this topic: Fang and Shao \cite{FS05, FS07}, Fang, Wang, and Wu \cite{FWW}, Gourcy and Wu \cite{GW}, Wang \cite{W02, W08}, and Wu and Zhang \cite{WZ}. The QTCI for diffusion-like semimartingales and reflected processes with respect to the uniform distance appeared in Pal \cite{pa} and was generalized by \"Ust\"unel \cite{U}. We refer the reader to \cite{pa} for more details on the existing literature. Similar ideas have been also put forward for discrete Markov chains. See, for example, the articles by Marton \cite{M2,M3}.

\subsection{RBM and the Skorokhod map} 

The reflected Brownian motions (RBMs) we are interested in have a constant drift vector, a constant diffusion matrix and are reflected whenever they hit the boundary of the positive orthant. On each face of the boundary of the latter the direction of reflection is constant. The theory of such processes is well-developed. In particular, their existence and pathwise uniqueness follows from the existence of a deterministic transformation mapping Brownian paths to the corresponding reflected paths. This is the so-called Skorokhod map, whose one-dimensional version is due to Skorokhod. We lift the following description from the article by Harrison and Reiman \cite{hr}. For more details and generalizations see \cite{DI}, \cite{DR1}, and \cite{DR2}.

Define $C([0,\infty), \rr^{K-1})$ as the space of continuous functions $x: [0, \infty) \rightarrow \rr^{K-1}$, endowed with the topology of the locally uniform convergence. For each such function $x$ we denote its component functions by $x_i$ for $i=1,\ldots, K-1$. Let $C_S$ denote the subset of functions in $C([0,\infty),\rr^{K-1})$ for which $x_i(0)\ge 0$, $i=1,\ldots,K-1$, and let $Q=(q_{ij})$ be a $(K-1)\times(K-1)$ matrix with non-negative entries, zeros on the diagonal and spectral radius strictly less than one. 

\begin{theorem}[\cite{hr}]\label{hrskoro}
For each $x \in C_S$ there exists a unique pair of functions $y, z \in C([0,\infty), \rr^{K-1})$ satisfying
\begin{eqnarray}\label{skorokhod}
z_j(t) &=& x_j(t) + y_j(t) - \sum_{i=1}^{K-1} q_{ij} y_i(t), \qquad t\ge 0,\\
z_j(t) &\ge& 0, \qquad t \ge0
\end{eqnarray}
for all $j=1,\ldots,K-1$ and such that for every $i=1,\ldots,K-1$, the function $y_i$ increases only at those times $t$ for which $z_i(t)=0$.\\ 
Moreover, suppose that the matrix $Q$ satisfies 
\eq\label{onenorm}
\norm{Q}_{cs}:=\max_{j=1,\dots,K-1} \sum_{i=1}^{K-1} q_{ij} < 1. 
\en
Then $y$ is the unique function in $C([0,\infty), \rr^{K-1})$ that for all $t\geq0$ satisfies the vector equation $y(t) = \sup_{0\le s \le t}[ y(s) Q - x(s) ]_+$, and is given by the limit (in the locally uniform topology) of the following iterative scheme:
\begin{eqnarray}
y^{[0]}(t) &\equiv& 0,\\
y^{[i+1]}(t) &=& \sup_{0\le s \le t}\left( y^{[i]}(s)Q - x(s)\right)_+, \qquad t \ge 0.
\end{eqnarray}
Hereby, the supremum and the positive part are taken componentwise.
\end{theorem}

It will become apparent later that the matrix $Q$ of our choice will be given by 
\begin{eqnarray*}
Q^{(K-1)}=\left(\begin{array}{cccccc}
0 & \frac{1}{2} & 0 & \hdots & \hdots & 0 \\
\frac{1}{2} & 0 & \frac{1}{2} & 0 & \hdots & 0\\
0 & \ddots & \ddots & \ddots & \ddots & \vdots \\
\vdots & \ddots & \frac{1}{2} & 0 & \frac{1}{2} & 0 \\
0 & \hdots & 0 & \frac{1}{2} & 0 & \frac{1}{2}\\
0 & \hdots & \hdots & 0 & \frac{1}{2} & 0 
\end{array}\right).
\end{eqnarray*}

Since $Q^{(K-1)}$ is a finite, irreducible, substochastic matrix, it is immediate (by adding an absorbing point) that the spectral norm of $Q^{(K-1)}$ is strictly less than one. As was shown in \cite{hr}, this implies that the corresponding Skorokhod map, which transforms the input path $x$ into the reflected path $z$ (or, the path of the ``local time'' $y$), is Lipschitz with respect to a variant of the norm $\|.\|_{T,2}$ defined in \eqref{T2norm}. We would like to obtain an explicit upper bound on its Lipschitz constant with respect to the norm $\|.\|_{T,2}$. 

Note that $Q^{(K-1)}$ does not satisfy the Assumption \eqref{onenorm} and, hence, the iterative scheme of Theorem \ref{hrskoro} cannot be applied directly to construct the corresponding Skorokhod map. The way to get around this is to define a $(K-1)\times(K-1)$ diagonal matrix
\begin{eqnarray*}
D^{(K-1)}=\left(\begin{array}{cccc} 
w^{(K-1)}\Big(\frac{1}{K}\Big) & 0 & \dots & 0 \\
0 & w^{(K-1)}\Big(\frac{2}{K}\Big) & \vdots & \vdots \\
\vdots & \hdots & \ddots & 0 \\  
0 & \dots & 0 & w^{(K-1)}\Big(\frac{K-1}{K}\Big)
\end{array}\right)
\end{eqnarray*}
with a strictly concave function $w^{(K-1)}: [0,1]\rightarrow[0,\infty)$ such that $w^{(K-1)}(u)=0$ if and only if $u\in\{0,1\}$. Viewing the set of equations in \eqref{skorokhod} for different values of $j$ as an equation between row vectors, and multiplying both sides of it by the matrix $\left[ D^{(K-1)} \right]^{-1}$ from the right, we obtain the equations
\begin{equation}\label{rescaling}
z'_j(t) = x'_j(t) + y'_j(t) - \sum_{i=1}^{K-1} r_{ij} y'_i(t), \qquad t\geq0
\end{equation}
for $j=1,\dots,K-1$ where $x'=x \left[ D^{(K-1)} \right]^{-1}$, $y'=y \left[ D^{(K-1)} \right]^{-1}$, $z'=z \left[ D^{(K-1)} \right]^{-1}$ and
\begin{equation}
R^{(K-1)} = (r_{ij})= D^{(K-1)} Q^{(K-1)} \left[ D^{(K-1)} \right]^{-1}.
\end{equation}

Note that the coordinates of $x'$, $y'$ and $z'$ can be computed by a simple rescaling of the corresponding coordinates of $x$, $y$ and $z$,
respectively. Moreover, it holds 
\begin{eqnarray}\label{csnormineq}
\|R^{(K-1)}\|_{cs}=\max_{l=1,\dots,K-1} \frac{w^{(K-1)}\Big(\frac{l-1}{K}\Big)+w^{(K-1)}\Big(\frac{l+1}{K}\Big)}{2w^{(K-1)}\Big(\frac{l}{K}\Big)}<1
\end{eqnarray}
due to the strict concavity of the function $w^{(K-1)}$. Thus, the Skorokhod map corresponding to the matrix $R^{(K-1)}$ can be obtained using the iterative scheme of Theorem \ref{hrskoro}. 

\bigskip

From now on we fix a terminal time $T>0$ and equip the space $C([0,T],\rr^{K-1})$ of continuous $\rr^{K-1}$-valued functions on $[0,T]$ with the norm
\begin{eqnarray}\label{sup2norm}
\|x\|_{T,2}:=\left(\frac{1}{K-1}\sum_{i=1}^{K-1} \sup_{0\leq t\leq T} x_i(t)^2\right)^{\frac{1}{2}}
\end{eqnarray}
where $x_i$, $i=1,\dots,K-1$ are the component functions of $x$. We write $\Phi^{Q^(K-1)}_L$ and $\Phi^{Q^{(K-1)}}_R$ for the maps that take a path $x\in C([0,T],\rr^{K-1})$ to the local time path $y\in C([0,T],\rr^{K-1})$ and the reflected path $z\in C([0,T],\rr^{K-1})$, respectively, corresponding to the Skorokhod problem with reflection matrix $Q^{(K-1)}$ defined above. The following lemma is one of the crucial steps in the proofs of Theorems \ref{slopetheorem} and \ref{thm_infinite}. It provides an upper bound on the Lipschitz constant of the Skorokhod map and is a significant improvement on an earlier attempt by Pal \cite{pa}.

{\lemma\label{lemma_lipschitz} For all natural numbers $K\geq2$ the map $\Phi^{Q^{(K-1)}}_L$ is Lipschitz on $C([0,T],\rr^{K-1})$ with respect to the norm $\|.\|_{T,2}$ defined in \eqref{sup2norm}. Moreover, its Lipschitz constant $Lip^{(K-1)}_L$ satisfies 
\begin{eqnarray}\label{lipineq}
Lip^{(K-1)}_L\leq 2\cdot(K-1)^{\frac{5}{2}}.
\end{eqnarray}}

\noindent{\it Proof.} 1) We fix a natural number $K$ as in the statement of the lemma and will prove the inequality (\ref{lipineq}) for that value of $K$. From the considerations preceeding the lemma we know that for each $g\in C([0,T],\rr^{K-1})$ the map $\Phi^{R^{(K-1)}}_L$ corresponding to the Skorokhod problem with reflection matrix $R^{(K-1)}$, evaluated at $g$, is given by the limit of the iterative scheme of Theorem \ref{hrskoro} with input $g$, restricted to the interval $[0,T]$. Now, let $\widetilde{g}$ be another function in $C([0,T],\rr^{K-1})$, for each $k\in\nn$ define $f^k$ and $\widetilde{f}^k$ as the results of the $k$-th step of the iterative scheme of Theorem \ref{hrskoro} with inputs $g$ and $\widetilde{g}$, respectively, and set $f$ and $\widetilde{f}$ for $\Phi^{R^{(K-1)}}_L(g)$ and $\Phi^{R^{(K-1)}}_L(\widetilde{g})$, respectively. Finally, define the norm $\|.\|_{T,max}$ on $C([0,T],\rr^{K-1})$ by
\begin{eqnarray}\label{normTmax}
\|x\|_{T,max}=\max_{i=1,\dots,K-1}\sup_{0\leq t\leq T} |x_i(t)|,
\end{eqnarray}
where $x_i$, $i=1,\dots,K-1$ denote the component functions of a function $x\in C([0,T],\rr^{K-1})$ as before. By applying the triangle inequality, the fact that the operation of taking the positive part is $1$-Lipschitz and the definition of the norm $\|.\|_{cs}$ in \eqref{onenorm}, one obtains the following chain of inequalities:
\begin{eqnarray*}
&&\|f^{k+1}-\widetilde{f}^{k+1}\|_{T,max}\\
&&=\max_{i=1,\dots,K-1}\sup_{0\leq t\leq T}\Big|\sup_{0\leq s\leq t}((f^k(s)R^{(K-1)})_i-g_i(s))_+ 
- \sup_{0\leq s\leq t}((\widetilde{f}^k(s)R^{(K-1)})_i-\widetilde{g}_i(s))_+\Big|\\
&&\leq\max_{i=1,\dots,K-1}\sup_{0\leq t\leq T}\sup_{0\leq s\leq t}\Big|((f^k(s)R^{(K-1)})_i-g_i(s))_+
-((\widetilde{f}^k(s)R^{(K-1)})_i-\widetilde{g}_i(s))_+\Big|\\
&&\leq\max_{i=1,\dots,K-1}\Big(\sup_{0\leq t\leq T} |(f^k(t)R^{(K-1)})_i-(\widetilde{f}^k(t)R^{(K-1)})_i|
+\sup_{0\leq t\leq T} |g_i(t)-\widetilde{g}_i(t)|\Big)\\
&&\leq\|R^{(K-1)}\|_{cs}\cdot\max_{i=1,\dots,K-1}\sup_{0\leq t\leq T} |f^k_i(t)-\widetilde{f}^k_i(t)|+\|g-\widetilde{g}\|_{T,max}.
\end{eqnarray*} 
Taking the limit $k\rightarrow\infty$ and rearranging terms we conclude
\begin{eqnarray}
\|f-\widetilde{f}\|_{T,max}\leq \frac{1}{1-\|R^{(K-1)}\|_{cs}}\cdot\|g-\widetilde{g}\|_{T,max}.
\end{eqnarray}

Recalling that the map $\Phi^{Q^{(K-1)}}_L$ can be obtained from the map $\Phi^{R^{(K-1)}}_L$ by a rescaling of the coordinates of $\rr^{K-1}$ according to the matrix $D^{(K-1)}$ (see equation \eqref{rescaling}), we deduce
\begin{eqnarray}\label{maxlip}
\|\Phi^{Q^{(K-1)}}_L(g)-\Phi^{Q^{(K-1)}}_L(\widetilde{g})\|_{T,max}
\leq\frac{\rho}{1-\|R^{(K-1)}\|_{cs}}\cdot\|g-\widetilde{g}\|_{T,max},
\end{eqnarray}  
where 
\begin{equation}
\rho=\frac{\max_{l=1,\dots,K-1} w^{(K-1)}(\frac{l}{K})}{\min_{l=1,\dots,K-1} w^{(K-1)}(\frac{l}{K})}. 
\end{equation}
Thus, the obvious equivalence of norms inequalities between the norms $\|.\|_{T,\max}$ and $\|.\|_{T,2}$ on $C([0,T],\rr^{K-1})$ show
\begin{eqnarray*}
&&\|\Phi^{Q^{(K-1)}}_L(g)-\Phi^{Q^{(K-1)}}_L(\widetilde{g})\|_{T,2}
\leq \|\Phi^{Q^{(K-1)}}_L(g)-\Phi^{Q^{(K-1)}}_L(\widetilde{g})\|_{T,max}\\
&&\leq \frac{\rho}{1-\|R^{(K-1)}\|_{cs}}\cdot\|g-\widetilde{g}\|_{T,max}
\leq \frac{(K-1)^{\frac{1}{2}}\rho}{1-\|R^{(K-1)}\|_{cs}}\cdot\|g-\widetilde{g}\|_{T,2}. 
\end{eqnarray*}

Since the functions $g$ and $\widetilde{g}$ were chosen arbitrarily in $C([0,T],\rr^{K-1})$, we conclude  
\begin{eqnarray}\label{lipbound}
Lip^{(K-1)}_L\leq\frac{(K-1)^{\frac{1}{2}}\rho}{1-\|R^{(K-1)}\|_{cs}}.
\end{eqnarray}
2) For $K=2$ or $K=3$ one can take $D^{(K-1)}$ to be the identity matrix of appropriate dimension, so that $\rho=1$ and $\|R^{(K-1)}\|_{cs}=\|Q^{(K-1)}\|_{cs}$ which immediately yields the inequality (\ref{lipineq}). To finish the proof for the case $K\geq4$, we need to bound the right-hand side in (\ref{lipbound}) from above for a suitable function $w^{(K-1)}$ in the definition of the matrix $D^{(K-1)}$. From \eqref{csnormineq} it is clear that
\begin{eqnarray}
\|R^{(K-1)}\|_{cs}\leq\max_{0\leq x\leq 1-\frac{2}{K}} \frac{w^{(K-1)}(x)+w^{(K-1)}\Big(x+\frac{2}{K}\Big)}{2w^{(K-1)}\Big(x+\frac{1}{K}\Big)}.
\end{eqnarray}

To define $w^{(K-1)}$, we choose a continuous function $v^{(K-1)}$ on $[0,1]$ such that $v^{(K-1)}(\frac{1}{2})=0$, $v^{(K-1)}(x)=-v^{(K-1)}(1-x)$, $x\in[0,1]$ and $v^{(K-1)}$ is strictly decreasing and convex on $[0,\frac{1}{2}]$. Then we set $w^{(K-1)}(x)=\int_0^x v^{(K-1)}(y)\;dy$ for all $x\in[0,1]$. In this case the convexity properties of $v^{(K-1)}$ and $w^{(K-1)}$ imply that
\begin{eqnarray}
\frac{d}{dx}\;\frac{w^{(K-1)}(x)+w^{(K-1)}\Big(x+\frac{2}{K}\Big)}{2w^{(K-1)}\Big(x+\frac{1}{K}\Big)}>0 
\end{eqnarray}
for $x\in[0,\frac{1}{2}-\frac{2}{K}]$ and that the opposite inequality holds for $x\in[\frac{1}{2},1-\frac{2}{K}]$. Moreover, if the latter derivative was zero at a point $x\in[\frac{1}{2}-\frac{2}{K},\frac{1}{2}-\frac{1}{K})$, then the strict concavity of the function $w^{(K-1)}$ would imply $v^{(K-1)}\Big(x+\frac{2}{K}\Big)<2v^{(K-1)}\Big(x+\frac{1}{K}\Big)-v^{(K-1)}(x)$. However, this is a contradiction to the fact that the graph of the function $v^{(K-1)}$ on $[\frac{1}{2},\frac{1}{2}+\frac{1}{K}]$ is never below the line connecting the points $\Big(\frac{1}{2},0\Big)$ and $\Big(\frac{1}{2}+\frac{1}{K},v^{(K-1)}\Big(\frac{1}{2}+\frac{1}{K}\Big)\Big)$. 

All in all, we conclude that the maximum in
\begin{eqnarray*}
\|R^{(K-1)}\|_{cs}=\max_{l=1,\dots,K-1} \frac{w^{(K-1)}\Big(\frac{l-1}{K}\Big)+w^{(K-1)}\Big(\frac{l+1}{K}\Big)}{2w^{(K-1)}\Big(\frac{l}{K}\Big)}
\end{eqnarray*}
is achieved when $\frac{l}{K}$ takes the value closest to $\frac{1}{2}$. If $K$ is even, the value $\frac{1}{2}$ is attained and, thus, the value of the latter maximum is given by
\begin{eqnarray*}
&&\frac{w^{(K-1)}\Big(\frac{1}{2}-\frac{1}{K}\Big)+w^{(K-1)}\Big(\frac{1}{2}+\frac{1}{K}\Big)}{2w^{(K-1)}\Big(\frac{1}{2}\Big)}
=\frac{\int_0^{\frac{1}{2}-\frac{1}{K}} v^{(K-1)}(y)\;dy+\int_0^{\frac{1}{2}+\frac{1}{K}} v^{(K-1)}(y)\;dy}
{2\int_0^{\frac{1}{2}} v^{(K-1)}(y)\;dy}\\
&&=1-\frac{\int_{\frac{1}{2}-\frac{1}{K}}^{\frac{1}{2}} v^{(K-1)}(y)\;dy}{\int_0^{\frac{1}{2}} v^{(K-1)}(y)\;dy}. 
\end{eqnarray*} 
If $K$ is odd, the maximum is achieved at $\frac{l}{K}=\frac{K-1}{2K}$ and for its value one has
\begin{eqnarray*}
&&\frac{w^{(K-1)}\Big(\frac{K-3}{2K}\Big)+w^{(K-1)}\Big(\frac{K+1}{2K}\Big)}{2w^{(K-1)}\Big(\frac{K-1}{2K}\Big)}
=1-\frac{\int_{\frac{K-3}{2K}}^{\frac{K-1}{2K}} v^{(K-1)}(y)\;dy-\int_{\frac{K-1}{2K}}^{\frac{K+1}{2K}} v^{(K-1)}(y)\;dy}
{2\int_0^{\frac{K-1}{2K}} v^{(K-1)}(y)\;dy}\\
&&\leq 1-\frac{\int_{\frac{K-2}{2K}}^{\frac{1}{2}} v^{(K-1)}(y)\;dy}{2\int_0^{\frac{1}{2}} v^{(K-1)}(y)\;dy}. 
\end{eqnarray*} 
Hereby, the last inequality is a consequence of the monotonicity and the symmetry properties of the function $v^{(K-1)}$, which show that we have decreased the numerator and increased the denominator of the fraction. The latter two calculations show that
\begin{eqnarray}\label{firstbnd}
\frac{\rho}{1-\|R^{(K-1)}\|_{cs}}\leq\frac{2\cdot\Big(\int_0^{\frac{1}{2}} v^{(K-1)}(y)\;dy\Big)^2}
{\Big(\int_{\frac{1}{2}-\frac{1}{K}}^{\frac{1}{2}}v^{(K-1)}(y)\;dy\Big)\cdot\Big(\int_0^{\frac{1}{K}} v^{(K-1)}(y)\;dy\Big)}.
\end{eqnarray}
For the unique function $v^{(K-1)}$ of the type described above for which $v^{(K-1)}(0)=K$ and $v^{(K-1)}(\frac{1}{K})=2$, and which is affine on $[0,\frac{1}{K}]$ and $[\frac{1}{K},\frac{1}{2}]$, the latter expression computes to
\begin{eqnarray}\label{secondbnd}
\frac{2}{\frac{2}{(K-2)K}\cdot\frac{K+2}{2K}}\leq 2\cdot(K-1)^2.
\end{eqnarray}
We remark at this point that the order of magnitude of the bound in terms of $K$ is optimal for this choice of the function $v^{(K-1)}$, since the non-negativity and the convexity of $v^{(K-1)}$ on $[0,\frac{1}{2}]$ together with $v^{(K-1)}(\frac{1}{2})=0$ show
\begin{eqnarray*}
\int_0^{\frac{1}{K}} v^{(K-1)}(y)\;dy\leq \int_0^{\frac{1}{2}} v^{(K-1)}(y)\;dy,\quad
\int_{\frac{1}{2}-\frac{1}{K}}^{\frac{1}{2}}v^{(K-1)}(y)\;dy\leq\Big(\frac{2}{K}\Big)^2\cdot \int_0^{\frac{1}{2}} v^{(K-1)}(y)\;dy.
\end{eqnarray*} 
Finally, combining the inequalities \eqref{lipbound}, \eqref{firstbnd} and \eqref{secondbnd}, we end up with the statement of the lemma. \ep 

\section{Concentration of the shape of the market weights}\label{sec_slope}

Consider the particle system in \eqref{ranksde} with $I=\{1,\dots,K\}$ for some $K\in\nn$ and a choice of constants $\delta_1 ,\dots,\delta_K$ that satisfy the following condition. Setting $\bar\delta=\frac{1}{K}\sum_{i=1}^K \delta_i$, one has
\eq\label{cond_delta}
\alpha_j:= \sum_{i=1}^j (\bar \delta - \delta_{K-i+1}) > 0,\quad j=1,\ldots,K-1.
\en  
In our analysis of the market weights we will use the following result from \cite{pp}.

\begin{lemma}\label{invatlas}
Under the condition \eqref{cond_delta} the process of spacings 
\eq\label{spacings}
(\xi_i(t):\;i=1,\dots,K-1):=(X_{(K-i+1)}(t)-X_{(K-i)}(t):\;i=1,\ldots,K-1),\quad t\geq0
\en
has a unique stationary distribution which is that of independent exponential random variables with rates $2\alpha_i$, for $i=1,\dots,K-1$. Moreover, the system of spacings is reversible at equilibrium. 
\end{lemma}

A situation in which the above lemma applies is given by the \textit{Atlas model} with $\delta_1=\delta>0$ and $\delta_i=0$ for all $i=2,\dots,K$. In that case one easily computes
\begin{equation}
\alpha_j= \sum_{i=1}^j \left( \frac{\delta}{K}  \right)= \frac{\delta j}{K},\quad j=1,\dots,K-1.
\end{equation}
Thus, under the stationary distribution we have 
\begin{eqnarray*}
&&\ev[\xi_j(t)]= \frac{K}{2\delta j},\quad j=1,\dots,K-1,\\ 
&&\ev[X_{(K-i+1)}(t)-X_{(K-j+1)}(t)]=\frac{K}{2\delta}\sum_{l=i}^{j-1} \frac{1}{l} \approx \frac{K}{2\delta} \log (j/i),\quad 1\leq i<j\leq K
\end{eqnarray*}
for all $t\geq0$.

\bigskip

Next, consider the linear regression in \eqref{linreg}. The corresponding ordinary least squares estimator for the slope parameter is given by the formula
\eq\label{whatisalphat}
\begin{split}
\alpha(t) &= \frac{\sum_{i=2}^J (\log i) \left( \log \mu_1(t) - \log \mu_i(t) \right)}{\sum_{i=2}^J \log^2i}=\frac{\sum_{i=2}^J (\log i) (X_{(K)}(t) - X_{(K-i+1)}(t))}{\sum_{i=2}^J \log^2 i}\\
&=\frac{\sum_{i=2}^J (\log i) \sum_{j=1}^{i-1} \xi_j(t)}{\sum_{i=2}^J \log^2 i}= \frac{\sum_{j=1}^{J-1} \xi_j(t) \sum_{i=j+1}^J \log i }{\sum_{i=2}^J \log^2 i}\\
&=  \frac{\sum_{i=1}^{J-1} (\log J!/i!) \xi_i(t)}{\sum_{i=1}^{J-1} \log^2(i+1) }.
\end{split}
\en
We use this formula in the proof below.
\bigskip

\noindent\textit{Proof of Theorem \ref{slopetheorem}.} The proof is broken down into several steps.

\begin{enumerate}
\item[\textit{Step 1.}] A Quadratic Transportation Cost Inequality for the process of spacings is derived in the Atlas model with $K$ particles.
\item[\textit{Step 2.}] Next, we assume $K$ to be very large compared to $J$. We prove a localization lemma that shows that the process $\alpha(t)$, $0\leq t\leq \delta^{-2}K$ is determined only by the particles corresponding to the $5J$ topmost indices with very high probability. 
\item[\textit{Step 3.}] Finally, we show that the law of the $5J$ particles in Step 2 is approximately that of another rank-based process, so that we can use the estimates obtained in step 1, with $K$ replaced by $5J$, to bound concentration of measure probabilities under the event constructed in step 2. 
\end{enumerate}

\bigskip

\noindent \textit{Step 1.} Consider the process $(X_{(1)}(t),\dots,X_{(K)}(t))$, $t\in[0,T]$ of ordered particles in the system \eqref{ranksde} with $I=\{1,\dots,K\}$. From Lemma 4 in \cite{pp} we know that there exist i.i.d. standard Brownian motions $\beta_1,\dots,\beta_K$ such that for all $i=1,\dots,K$ it holds
\begin{eqnarray}  
dX_{(i)}(t)=\delta_i\;dt+d\beta_i(t)+\frac{1}{2}dL_{(i-1,i)}(t)-\frac{1}{2}dL_{(i,i+1)}(t),\quad t\in[0,T]
\end{eqnarray}
with $L_{i,i+1}(t)$, $t\in[0,T]$ being the local time process at zero of the process $X_{(i+1)}(t)-X_{(i)}(t)$, $t\in[0,T]$ for $i=1,\dots,K-1$ and the convention $L_{0,1}(t)=L_{K,K+1}(t)=0$ for all $t\in[0,T]$. 

Hence, for all $i=1,\dots,K-1$ one has the dynamics   
\begin{eqnarray*}
d(X_{(i+1)}(t)-X_{(i)}(t))&=&(\delta_{i+1}-\delta_i)dt+d\beta_{i+1}(t)-d\beta_i(t)\\
&+&dL_{(i,i+1)}(t)-\frac{1}{2}dL_{(i+1,i+2)}(t)-\frac{1}{2}dL_{(i-1,i)}(t)
\end{eqnarray*}   
on $[0,T]$. In particular, we can conclude from this representation as in section 2 of \cite{pp} that the process 
\begin{eqnarray*}
(X_{(2)}(t)-X_{(1)}(t),\dots,X_{(K)}(t)-X_{(K-1)}(t)),\quad t\in[0,T]
\end{eqnarray*}
is a reflected Brownian motion in the $(K-1)$-dimensional positive orthant with reflection matrix $Q^{(K-1)}$ in the sense of section 1 in \cite{hr}. By reversing the labeling we see that the process $\xi(t):=(\xi_1(t), \dots, \xi_{K-1}(t))$, $t\in[0,T]$ is also an RBM in the positive orthant with reflection matrix $Q^{(K-1)}$. Thus, the process $\xi$ can be obtained as the image of the process $\gamma^*(t):=(\beta_{i+1}(t) + \delta_{i+1} t - \beta_i(t) - \delta_i t, \; i=1,\dots, K-1)$, $t\in[0,T]$ under the map $\Phi^{Q^{(K-1)}}_R$.  

By Theorem 6 in \cite{pa} the process of independent Brownian motions $(\beta_1(t) + \delta_1 t, \dots, \beta_K(t) + \delta_K t )$, $t\in[0,T]$ satisfies a QTCI with respect to the norm $\|.\|_{T,2}$ with the constant $4K^{-1}T$. Moreover, the map that takes the vector of these Brownian motions to the process $\gamma^*(t)$, $t\in[0,T]$ is Lipschitz with respect to the norm $\norm{\cdot}_{T,2}$ with Lipschitz constant $2\sqrt{K/(K-1)}$. Hence, by Lemma \ref{lemmapushforward} the process $\gamma^*(t)$, $t\in[0,T]$ satisfies a QTCI with respect to the norm $\norm{\cdot}_{T,2}$ with the constant 
\eq\label{firstlip}
C^*_K(0):= 16 K^{-1} T \frac{K}{K-1} = 16 (K-1)^{-1} T.
\en

Now, we will use Lemma \ref{lemma_lipschitz}. Consider the equation \eqref{skorokhod} with the matrix $Q^{(K-1)}$ for a fixed value of $t$ and two different unconstrained processes $x$ and $\tilde x$. Writing $y,\tilde y,z, \tilde z$ for $\Phi^{Q^{(K-1)}}_L(x)$, $\Phi^{Q^{(K-1)}}_L(\tilde x)$, $\Phi^{Q^{(K-1)}}_R(x)$, $\Phi^{Q^{(K-1)}}_R(\tilde x)$, respectively, we easily deduce
\[
\abs{z_j(t) - \tilde z_j(t)} \le \abs{x_j(t) - \tilde x_j(t)} + \Big(\max_{i=1,\dots,K-1}\abs{y_i-\tilde y_i}\Big)\cdot \underbrace{\sum_{i=1}^{K-1} \abs{(I^{(K-1)}-Q^{(K-1)})_{ij}}}_{\leq 2},\quad t\in[0,T], 
\]
where $I^{(K-1)}$ is a $(K-1)\times(K-1)$ identity matrix. Hence,
\begin{equation}\label{Tmaxbound1}
\norm{z-\tilde z}_{T,\max} \le \norm{x-\tilde x}_{T,\max} + 2 \norm{y - \tilde y}_{T,\max}.
\end{equation}
On the other hand, combing the bounds \eqref{maxlip}, \eqref{firstbnd} and \eqref{secondbnd} we get
\begin{equation}\label{Tmaxbound2}
\norm{y - \tilde y}_{T,\max} \le 2(K-1)^2 \norm{x-\tilde x}_{T,\max}
\end{equation}
for all $K\geq 4$. By the same argument as in the beginning of step 2 in the proof of Lemma \ref{lemma_lipschitz} the same inequality is true for $K=2$ and $K=3$. Putting the inequalities \eqref{Tmaxbound1} and \eqref{Tmaxbound2} together we obtain  
\begin{equation}\label{Rlipbound}
\norm{z-\tilde z}_{T,\max} \le \left(1+4(K-1)^2\right)\cdot \norm{x - \tilde x}_{T,\max}\le 5(K-1)^2 \norm{x-\tilde x}_{T,\max}
\end{equation}
for all $K\geq2$. Recall that by applying the map $\Phi^{Q^{(K-1)}}_R$ to the paths of the process $\gamma^*$ and by reversing the order of the coordinates thereafter one gets the process $\xi$. Combining this observation with Lemma \ref{lemmapushforward}, \eqref{firstlip} and \eqref{Rlipbound}, we see that the process $\xi(t)$, $t\in[0,T]$ satisfies a QTCI with respect to the norm $\norm{\cdot}_{T,\max}$ with the constant 
\eq\label{lipxi}
C^*_K := 400(K-1)^3 T.
\en

Now, we restrict ourselves to the first $(J-1)$ coordinates of the process $\xi$, since only those appear in \eqref{whatisalphat}. By Lemma \ref{lemmapushforward} the vector-valued process $(\xi_1(t), \ldots, \xi_{J-1}(t))$, $t\in[0,T]$ also satisfies a QTCI with respect to the norm $\norm{\cdot}_{T,\max}$ with the constant $C_K^*$.
\bigskip

\noindent\textit{Step 2.} Consider the formula for $\alpha(t)$ in \eqref{whatisalphat}. The value of $\alpha(t)$ depends only on the top $J$ spacings, whereby we have assumed that $J$ is very much smaller than $K$. In this case, for a large enough $m$ and with high probability, the top $J$ processes during the time interval $[0,T]$ are identical to the top $J$ processes among the processes which start off at the top $J+m$ positions at time zero. The following lemma makes this idea precise. 

\begin{lemma}\label{localprob}
Consider the particle system of Theorem \ref{slopetheorem} and for all $m\in\nn$ define $\sigma_m$ as the first time $t$ at which, for some $i \ge J+m$ and some $1\le j \le J$, it holds $X_{K-i+1}(t)= X_{(K-j+1)}(t)$. Then, for all $K \ge 30$, we have
\eq\label{sigmaestimate}
\pp\left( \sigma_{4J+1} \le \delta^{-2}K  \right) \le 6 J^2 K^{-3/2} e^{-K/50}.
\en
\end{lemma}

\noindent\textit{Proof of the Lemma.} For any $T>0$, the event $\{ \sigma_m \le T \}$ implies that for some $i\ge J+m$ and some $1\le j\le J$, the processes $X_{K-i+1}$ and $X_{K-j+1}$ cross paths during the time interval $[0,T]$. Using the union bound and bounding the drift of the lower particle by the constant $\delta$, we have the following estimate:
\eq\label{probsigma}
\begin{split}
\pp\left( \sigma_{m} \le T  \right) &\le \sum_{i=J+m}^K \sum_{j=1}^J \;\pp\left(  \sup_{0\le t \le T} \left( W_{K-i+1}(t)-W_{K-j+1}(t) \right) \ge -\delta T + X_{K-j+1}(0) - X_{K-i+1}(0) \right)\\
&\le J \sum_{i=J+m}^K \pp\left(  \sup_{0\le t \le T} \left( W_{K-i+1}(t) - W_{K-J+1}(t) \right) \ge -\delta T + X_{K-J+1}(0) - X_{K-i+1}(0) \right).
\end{split}
\en

Next, we use the fact that the supremum of a standard Brownian motion up to time $T$ has the same law as the absolute value of a normal random variable with mean zero and variance $T$. Thus, with a standard normal random variable $Z$ one has 
\[
\begin{split}
\pp&\left(  \sup_{0\le t \le T} \left( W_{K-i+1}(t) - W_{K-J+1}(t) \right) \ge -\delta T + X_{K-J+1}(0) - X_{K-i+1}(0) \right)\\
&= \pp\left( \sqrt{2T}\abs{Z}    \ge \frac{K}{\delta}\log(i/J) - \delta T\right)
=2\phibar\left(  \left(\frac{K}{\delta\sqrt{2T}}\log(i/J) - \delta \sqrt{\frac{T}{2}}\right)_+ \right).  
\end{split}
\]

Here, $\phibar$ is the one minus the cumulative distribution function of a standard normal random variable.

Plugging this into \eqref{probsigma}, we get
\eq\label{sigmabnd2}
\begin{split}
\pp\left( \sigma_{m} \le T  \right) &\le 2J \sum_{i=J+m}^K \phibar\left(  \left(\frac{K}{\delta\sqrt{2T}}\log(i/J) - \delta \sqrt{\frac{T}{2}} \right)_+\right)\\
&\le 2J \int_{J+m-1}^\infty \phibar\left(  \left( \frac{K}{\delta\sqrt{2T}}\log(x/J) - \delta \sqrt{\frac{T}{2}} \right)_+\right) dx.
\end{split}
\en
We note that for $m\geq 4J+1$ and $T=\delta^{-2}K$ we may omit taking the positive part in the latter formula. 

Now, recall the well-known inequalities (see \cite[p.~298]{AS})
\eq\label{gausstail}
\frac{2\phi(y)}{y+ \sqrt{y^2+4}}\le \phibar(y) \le \frac{1}{y} \phi(y),\quad y>0, 
\en
where $\phi$ is the density function of a standard normal random variable. From the latter inequality it follows that the function $\log \phibar(y)$ is concave on $(0,\infty)$. To see this, note that
\[
\left(\log \phibar(y)\right)'' = \left(  y - \frac{\phi(y)}{\phibar(y)}   \right) \frac{\phi(y)}{\phibar(y)} < 0,\quad y>0.
\]
In other words, for any positive real numbers $a$ and $b$, we have
\[
\phibar\left(  \frac{a+b}{2} \right) \ge \sqrt{\phibar(a) \phibar(b)}.
\]

Taking $T=\delta^{-2}K$ and $m\geq 4J+1$, and choosing
\[
a= \frac{K}{\delta\sqrt{2T}}\log(x/J) - \delta \sqrt{\frac{T}{2}} , \quad b= \delta \sqrt{\frac{T}{2}}
\]
for $x\geq J+m-1$ in the latter inequality, we get
\[
\phibar\left(\frac{K}{\delta\sqrt{2T}}\log(x/J) - \delta \sqrt{\frac{T}{2}} \right) 
\le \left[\phibar\left( \delta \sqrt{\frac{T}{2}}\right)\right]^{-1}\cdot \left[\phibar\left(  \frac{K}{2\delta\sqrt{2T}} \log(x/J) \right)\right]^2
\]
for all $x\geq J+m-1$. Substituting in \eqref{sigmabnd2}, we deduce
\[
\begin{split}
\pp\left( \sigma_{m} \le T  \right) &\le 2J \left[\phibar\left( \delta \sqrt{\frac{T}{2}}  \right)\right]^{-1} \int_{J+m-1}^\infty \left[\phibar\left(  \frac{K}{2\delta\sqrt{2T}} \log(x/J) \right)\right]^2 dx\\
&=\frac{4J^2\delta\sqrt{2T}}{K} \left[\phibar\left( \delta \sqrt{\frac{T}{2}}  \right)\right]^{-1} \int_{a^*}^\infty \exp\left(2\delta\sqrt{2T}y/K\right) \Big[\phibar(y)\Big]^2  dy,
\end{split}
\]
where $a^*=\log[(J+m-1)/J]\; K/(2\delta\sqrt{2T})$. Hence, using $\phibar(y)\leq y^{-1}\phi(y)$, $y>0$, we obtain
\[
\begin{split}
\pp\left( \sigma_{m} \le T  \right) &\le \frac{2J^2\delta \sqrt{2T}}{\pi Ka^{*2}} \left[\phibar\left( \delta \sqrt{\frac{T}{2}}  \right)\right]^{-1} \int_{a^*}^\infty \exp\left(2\delta\sqrt{2T}y/K - y^2\right)  dy\\
&\le \frac{2J^2\delta\sqrt{2T}}{\sqrt{\pi}K a^{*2}} \left[\phibar\left( \delta \sqrt{\frac{T}{2}}  \right)\right]^{-1} e^{2T\delta^2/K^2} \phibar\left( \sqrt{2}\left(  a^* - \frac{\delta\sqrt{2T}}{K} \right)   \right).
\end{split}
\]

At this point we take $m=4J+1$ and substitute $\delta^{-2}K$ for $T$. In this case, one has $a^*=\frac{\log 5}{2}\cdot\sqrt{\frac{K}{2}}$ and the latter bound simplifies to
\[
\pp\left( \sigma_{4J+1} \le \delta^{-2}K  \right) \le\frac{\widetilde{C}_1}{K^{3/2}}\left[\phibar\left(\sqrt{\frac{K}{2}}\right)\right]^{-1}
\phibar\left(\frac{\log 5}{2}\cdot\sqrt{K} - \frac{2}{\sqrt{K}}\right).
\]
Hereby, we have set 
\[
\widetilde{C}_1=\frac{16\sqrt{2}J^2}{\sqrt{\pi}(\log 5)^2} e^{2K^{-1}}.
\]

Using the bounds in \eqref{gausstail} we obtain
\[
\begin{split}
&\phibar\left( \sqrt{\frac{K}{2}}\right) 
\ge \phi\left( \sqrt{\frac{K}{2}}\right)\cdot\frac{2K^{-1/2}}{\left(1/\sqrt{2} + \sqrt{1/2+ 4 K^{-1}} \right)},\\
&\phibar\left(  \frac{\log 5}{2} \sqrt{K} - \frac{2}{\sqrt{K}} \right) 
\le \frac{\phi\left( \frac{\log 5}{2} \sqrt{K} - \frac{2}{\sqrt{K}}    \right)}{\left(  \frac{\log 5}{2} \sqrt{K} - \frac{2}{\sqrt{K}}   \right)}
= \frac{\phi\left(\widetilde{C}_2 \sqrt{K}\right)}{\widetilde{C}_2 \sqrt{K}}
\end{split}
\]
for all $K\geq 30$, where
\[
\widetilde{C}_2=\frac{\log 5}{2}-\frac{2}{K}.
\]

Hence, combining all our estimates, we end up with
\[
\pp\left( \sigma_{4J+1} \le \delta^{-2}K  \right) 
\le \frac{\widetilde{C}_1}{K^{3/2}}\cdot\frac{1/\sqrt{2} + \sqrt{1/2+ 4 K^{-1}}}{2\widetilde{C}_2} 
\exp\left(\frac{K}{4}  - \frac{\widetilde{C}_2^2}{2} K \right).
\]
To finish the proof of the lemma, it suffices to observe that the right-hand side of the latter inequality is bounded above by $6 J^2 K^{-3/2} e^{-K/50}$ for all $K \ge 30$. 

\bigskip

\noindent\textit{Step 3.} Consider the particle system of Theorem \ref{slopetheorem}. We proceed with another \textit{localization} lemma.

\begin{lemma}\label{local2}
Define $\tilde{\sigma}_m$ to be the first time $t$ at which, for some $i \ge m$, we have $X_i(t) = X_{(1)}(t)$. Then
\[
\pp\left(  \tilde{\sigma}_{2K/3+1} \le \delta^{-2}K  \right) \le 11\sqrt{K} e^{-K/500}.
\]  
\end{lemma}

\noindent\textit{Proof of the Lemma.} The proof is very similar to that of Lemma \ref{localprob} and we only outline the argument. With $m=2K/3+1$, and letting $Z$ be a standard normal random variable and $\overline{\Phi}$ be one minus its cumulative distribution function as before, we have
\[
\begin{split}
\pp\left( \tilde{\sigma}_m \le \delta^{-2}K  \right) &\le \sum_{i=m}^K \pp\left( \sup_{0\le t \le \delta^{-2}K}\left( W_1(t) - W_i(t) \right) \ge -\delta^{-1} K + X_i(0) - X_1(0)   \right)\\
&= \sum_{i=m}^K \pp\left( \delta^{-1} \sqrt{2K}\abs{Z} \ge \delta^{-1}K\log\left(\frac{K}{K-i+1}\right) -\delta^{-1} K \right)\\
&= \sum_{i=m}^K 2\phibar\left(\left(   \sqrt{\frac{K}{2}} \log\left(\frac{K}{K-i+1}\right) - \sqrt{\frac{K}{2}}   \right)_+\right)\\
&= 2\sum_{i=1}^{K-m+1} \phibar\left(   \sqrt{\frac{K}{2}} \log\left(\frac{K}{i}\right) - \sqrt{\frac{K}{2}}   \right).
\end{split}
\]
Note that we could drop the positive part in the last identity due to our assumption $m=2K/3+1$. Now, we use the fact that for $i \le K- (2K/3+1)+1$, one has $\log(K/i)-1\ge \log3-1\ge .09$, together with the second inequality in \eqref{gausstail} to obtain 
\[
\pp\left( \tilde{\sigma}_{2K/3+1} \le \delta^{-2}K  \right)\le \frac{2K}{3} \phibar\left( .09 \sqrt{\frac{K}{2}}  \right) \le \frac{2\sqrt{2K}}{.27} e^{-.0081 K/4}\le 11 \sqrt{K} e^{-K/500}.
\]

This finishes the proof of the Lemma. \ep

\bigskip

To complete the proof of Theorem \ref{slopetheorem}, we recall from \cite{pp} that a weak solution for the Atlas model as in Theorem \ref{slopetheorem} can be obtained by the following application of Girsanov's Theorem. Let $Z_1, \ldots, Z_K$ be independent Brownian motions such that $Z_i(0)=X_i(0)$, $i=1,\dots,K$. Set $T=\delta^{-2}K$ and let $Q^0$ denote their joint law during the time interval $[0,T]$ on the canonical sample space of continuous $\rr^K$-valued functions on $[0,T]$ with the usual Brownian filtration. Consider the martingale
\[
M(t) = \sum_{i=1}^K \int_0^t 1_{\{ Z_i(s) = Z_{(1)}(s)   \}} dZ_i(s),\qquad t\geq0
\]
Note that its quadratic variation at any fixed time $t\geq0$ is given by
\[
\iprod{M}(t) = \left(  \sum_{i=1}^K \int_0^t 1_{\{ Z_i(s) = Z_{(1)}(s)   \}} \right) ds = t,
\]
since two independent Brownian particles can simultaneously be the leftmost ones only on a set of Lebesgue measure zero. Hence, by L\'evy's characterization of Brownian motion, the process $M$ is a standard Brownian motion under $Q^0$. If we now change the measure to $Q^\delta$ according to the formula 
\eq\label{changeofmeasure}
\frac{d Q^\delta}{d Q^0}= \exp\left( \delta M(T) - \frac{\delta^2 T}{2}   \right),
\en
then, under $Q^\delta$, the law of the process $(Z_1(t),\ldots, Z_K(t))$, $t\in[0,T]$ is that of the Atlas model during the time interval $[0,T]$.

Let $\mathcal{F}_T$ denote the $\sigma$-algebra generated by the entire process $(Z_1(t), \ldots, Z_K(t)), \; t \in [0,T]$. Clearly, one has the decomposition
\[
\mathcal{F}_T = \mathcal{G}_T \vee \mathcal{H}_T,
\]
where $\mathcal{G}_T$ is the $\sigma$-algebra generated by the top $5J$ indexed coordinate processes
\[
(Z_{K-5J+1}(t), \ldots, Z_K(t)),\; t \in [0,T] ,
\]
$\mathcal{H}_T$ is the $\sigma$-algebra generated by the rest of the coordinate processes, and $\vee$ refers to the smallest $\sigma$-algebra containing the two.   

By our assumption, $5J < K/3 $. Hence, the process $(Z_1(t), \ldots, Z_{2K/3+1}(t)),\; t \in[0,T]$ is measurable with respect to $\mathcal{H}_T$. For any fixed $t\in[0,T]$ define 
\[
\newz(t) = \min_{i=1,\dots,2K/3+1} Z_i(t).
\]
Then, the process $\newz(t)$, $t\in[0,T]$ is also measurable with respect to $\mathcal{H}_T$.

Now, consider an arbitrary $\mathcal{G}_T$-measurable function $F$ such that $0\le F \le 1$. To simplify the notation, we will denote expectations with respect to the measures $Q$ and $Q^\delta$ by $Q(\cdot)$ and $Q^\delta(\cdot)$, respectively. By Lemma \ref{local2} and the change of measure formula \eqref{changeofmeasure} we have
\eq\label{measurechange1}
\begin{split}
Q^{\delta}\left( F \right) 
&\leq Q^{\delta}\left( F 1_{\{\tilde{\sigma}_{2K/3+1}>T\}}\right) + Q^{\delta}\left( 1_{\{\tilde{\sigma}_{2K/3+1}\leq T\}}\right)\\
&\le Q^0\left(  Fe^{\delta M(T) - \delta^2T/2}1_{\{ \tilde{\sigma}_{2K/3+1} > T \}} \right) + 11 \sqrt{K} e^{-K/500}.
\end{split}
\en

Now, on the set $\{ \tilde{\sigma}_{2K/3 +1} > T \}$, the process $M(t)$, $t\in[0,T]$ is identical to the process $\tilde{M}(t)$, $t\in[0,T]$, where the latter is defined by
\[
\tilde{M}(t) = \sum_{i=1}^{2K/3+1}\int_0^t 1_{\left\{ Z_i(s)= \newz(s)  \right\}} dZ_i(s), \quad t\in[0,T]. 
\]

Hence, it holds 
\[
Q^0\left(  Fe^{\delta M(T) - \delta^2T/2}1_{\{ \tilde{\sigma}_{2K/3+1} > T \}} \right) 
=Q^0\left(  Fe^{\delta\tilde{M}(T) - \delta^2T/2}1_{\{ \tilde{\sigma}_{2K/3+1} > T \}} \right) 
\le Q^0\left(  Fe^{\delta\tilde{M}(T) - \delta^2T/2} \right).
\]

Note that $\tilde{M}(T)$ is measurable with respect to $\mathcal{H}_T$ while $F$ is measurable with respect to $\mathcal{G}_T$.
Moreover, under $Q^0$, the $\sigma$-algebras $\mathcal{G}_T$ and $\mathcal{H}_T$ are independent of each other. Using this observation and $\iprod{\tilde{M}}(T)=T$, we obtain
\[
Q^0\left(  Fe^{\delta\tilde{M}(T) - \delta^2T/2} \right) = Q^0(F) Q^0\left( e^{\delta\tilde{M}(T) - \delta^2T/2}  \right)= Q^0(F).
\]

Combining this with the previous inequality and \eqref{measurechange1}, we get
\eq\label{measurechange2}
Q^{\delta}\left( F \right) \le Q^0(F) + 11 \sqrt{K} e^{-K/500}.
\en
For the rest of the argument we will assume that the particles indexed by the top $5J$ indices are independent Brownian motions starting from their respective initial conditions, the idea being that all probabilities under the actual measure $Q^\delta$ can be bounded from above as in \eqref{measurechange2}.

\bigskip

Back to the $K$ particles Atlas model, consider the event $\{ \sigma_{4J+1} > \delta^{-2}K \}$ as in Lemma \ref{localprob}. On this event, during the time interval $[0,\delta^{-2}K]$ the top $J$ processes are identical to the top $J$ processes among those that started at the top $5J$ positions at time zero. Let $Y_1',\dots,Y_{5J}'$ be the ranked processes $X_{K-5J+1},\dots,X_K$ in the increasing order. Also, set
\[
\alpha'(t) = \frac{\sum_{i=2}^{J}(\log i)(Y'_{5J}(t) - Y'_{5J+1-i}(t)) }{\sum_{i=2}^{J} \log^2 i},\quad t\geq0.
\]  
Then, by Lemma \ref{localprob},
\eq\label{alphaalphaprime}
\pp\left( \alpha(t) = \alpha'(t), \; \text{for all}\; 0\le t \le \delta^{-2}K   \right)  \ge 1 - 6 J^2 K^{-3/2} e^{-K/50}.
\en

We now prove a concentration of measure property for $\alpha'(t)$, $t\in[0,\delta^{-2}K]$. Relying on \eqref{measurechange2}, we can assume first that $X_{K-5J+1}, \ldots, X_K$ evolve according to independent standard Brownian motions. If we let 
\[
(\xi'_i(t):\;i=1,\dots,J-1)=(Y'_{5J-i+1}(t) - Y'_{5J-i}(t):\;i=1,\dots,J-1),\quad t\in[0,\delta^{-2}K],
\]
then $(\xi'_1(t),\dots,\xi'_{J-1}(t))$, $t\in[0,\delta^{-2}K]$ can be viewed as a vector of $(J-1)$ component processes of a reflected Brownian motion in the $(5J)$-dimensional positive orthant with zero drift vector and a constant diffusion matrix. 

By \eqref{lipxi} and the paragraph following it, we know that the process $\xi'(t)$, $t\in[0,\delta^{-2}K]$ satisfies a QTCI with respect to the norm $\|.\|_{\delta^{-2}K,max}$ with the constant $\mu' J^3\delta^{-2}K$. Here, we have abbreviated $5^3\times 400$ by $\mu'$. 

Note that for any fixed $t\in[0,\delta^{-2}K]$, $\alpha'(t)$ can be written in terms of $\xi'(t)$ as
\[
\alpha'(t) = \frac{\sum_{i=1}^{J-1} \xi_i'(t) \sum_{j=i+1}^J\log j}{\sum_{i=2}^J \log^2 i}
= \frac{\sum_{i=1}^{J-1} \xi_i'(t) \log(J!/ i!)}{\sum_{i=2}^J \log^2 i}.
\]

Hence, the function that takes the paths of $\xi'(t)$, $t\in[0,\delta^{-2}K]$ to the paths of $\alpha'(t)$, $t\in[0,\delta^{-2}K]$ is Lipschitz with respect to $\norm{\cdot}_{\delta^{-2}K,\max}$ norms with the Lipschitz constant
\[
\tilde{C}_{\alpha}(J) = \frac{\sum_{i=1}^{J-1} \log(J!/ i!)}{\sum_{i=2}^J \log^2 i}.
\]
It follows that the random variable $\overline{\alpha}':=\sup_{0\le t \le \delta^{-2}K} \alpha'(t)$ can be viewed as the image of $\xi'(t)$, $t\in[0,\delta^{-2}K]$ under a Lipschitz function with the Lipschitz constant $\tilde{C}_{\alpha}(J)$. Thus, by Lemma \ref{lemmapushforward} its law satisfies a QTCI with the constant
\[
\tilde{C}_{\alpha}(J)^2\mu' J^3 \delta^{-2}K.
\]

Let $m_{\alpha}$ be the median of $\overline{\alpha}'$ under $Q^0$. Setting $\mu=8\mu'$, we deduce from Theorem \ref{preq} and \eqref{measurechange2}:
\[
Q^\delta\left(  \overline{\alpha}' > m_{\alpha} + r\sqrt{K} \right)\le \exp\left( -\frac{r^2\delta^2}{\mu \tilde{C}_{\alpha}(J)^2 J^3}  \right) + 11\sqrt{K} e^{-K/500}
\]
for all $r$ greater than a constant depending only on $J$ and $\delta$. Combining this estimate with \eqref{alphaalphaprime} (and bounding $J$ by $K$) we get
\[
Q^\delta\left(  \overline{\alpha} > m_{\alpha} + r\sqrt{K} \right)\le \exp\left( -\frac{r^2\delta^2}{\mu \tilde{C}_{\alpha}(J)^2 J^3}  \right) + 11\sqrt{K} e^{-K/500} + 6\sqrt{K}e^{-K/50}.
\]
The observation that the sum of the last two summands is smaller than the first summand for all sufficiently large $K$ yields
\[
Q^\delta\left(  \overline{\alpha} > m_{\alpha} + r\sqrt{K} \right)\le 2 \exp\left( -\frac{r^2\delta^2}{\mu \tilde{C}_{\alpha}(J)^2 J^3}  \right).
\]

We note that $m_{\alpha}$ is not the median of $\overline\alpha$ under $Q^\delta$. However, by \eqref{measurechange2} and \eqref{alphaalphaprime} one has
\[
\begin{split}
Q^\delta\left(  \overline\alpha< m_{\alpha} \right) &\le Q^0\left(  \overline{\alpha} < m_{\alpha} \right)  + 11\sqrt{K} e^{-K/500} \\
&\le Q^0\left(  \overline\alpha' < m_{\alpha} \right)  + 11\sqrt{K} e^{-K/500} + 6\sqrt{K}e^{-K/50}\\
&\le 1/2 + 17\sqrt{K} e^{-K/500}.
\end{split}
\]
This completes the proof of Theorem \ref{slopetheorem}. \ep

\section{The infinite rank-based system}\label{sec_infinite}

This section is devoted to the proof of Theorem \ref{thm_infinite}. The first step in the proof is to understand the dynamics of the process $(X_{(1)}(t),\dots,X_{(n)}(t))$, $t\in[0,T]$ of the $n$ leftmost particles in the particle system of Theorem \ref{thm_infinite}.

{\lemma\label{lemmastopping} There exist stopping times $0=\tau_N\leq\tau_{N+1}\leq\dots$ such that the following is true. 
\begin{enumerate}[(a)]
\item $\lim_{m\rightarrow\infty}\tau_{N+m}=\infty$ with probability one.
\item For each $m\in\nn$ there exists a system of i.i.d. standard Brownian motions $\beta^{(m)}_1,\dots,\beta^{(m)}_n$ such that for all $i=1,\dots,n$ one has the dynamics
\begin{eqnarray*}
dX_{(i)}(t\wedge\tau_{N+m})=1_{\{\tau_{N+m}\geq t\}}\delta_i\;dt+d\beta^{(m)}_i(t\wedge\tau_{N+m})+\frac{1}{2}dL_{(i-1,i)}(t\wedge\tau_{N+m})\\
-\frac{1}{2}dL_{(i,i+1)}(t\wedge\tau_{N+m})
\end{eqnarray*}
on $[0,T]$. Hereby, $a\wedge b$ denotes $\min(a,b)$ for any two real numbers $a,b$.  
\end{enumerate}}

\noindent {\it Proof.} We define inductively the sets $\Lambda_N\subset \Lambda_{N+1}\subset\dots$ and the stopping times $0=\tau_N\leq\tau_{N+1}\leq\dots$ by
\begin{eqnarray}
\Lambda_{N+m}&=&\{k\geq1|\;\exists\; 1\leq i\leq N,\;0\leq s\leq\tau_{N+m}:\;X_k(s)=X_{(i)}(s)\},\\
\tau_{N+m+1}&=&\inf\{s\geq\tau_{N+m}|\;\exists\; 1\leq i\leq N,\;k\notin \Lambda_{N+m}:\;X_k(s)=X_{(i)}(s)\}
\end{eqnarray}
for all $m=0,1,\dots$. 

The proof of Proposition 3.1 in \cite{sh} shows that, with probability one, it holds $\lim_{m\rightarrow\infty}\tau_{N+m}=\infty$ and the sets $\Lambda_{N+m}$ are finite for all $m=0,1,\dots$. Moreover, the same proof implies that for each such number $m$ the paths of the process  $X_{(1)}(t\wedge\tau_{N+m}),\dots,X_{(n)}(t\wedge\tau_{N+m})$, $t\in[0,T]$ are given by the paths of the $n$ leftmost particles in a particle system as in \eqref{ranksde} with $I=\{1,\dots,|\Lambda_{N+m}|\}$, which is stopped at time $\tau_{N+m}$. Hence, by Lemma 4 in \cite{pp} we conclude that assertion (b) of the lemma is true for our choice of the stopping times $0=\tau_N\leq\tau_{N+1}\leq\dots$. \ep

\bigskip

Next, fix a $K\in\nn$ and let $X'_{(1)},\dots,X'_{(K)}$ be the ranked particles in the system \eqref{ranksde} with $I=\{1,\dots,K\}$. Also, let $L'_{(1,2)},\dots,L'_{(K-1,K)}$ be the local time processes at zero of the spacings processes in that system. Recall the definition of the norm $\|.\|_{T,2}$ in \eqref{T2norm}. From Lemma \ref{lemma_lipschitz} we can deduce the following concentration of measure property of the finite particle system.    
  
{\cor\label{corfinite} Let $A$ and $B$ be measurable subsets of $(C([0,T],\rr^{n-1}),\|.\|_{T,2})$ such that
\begin{eqnarray}
&&\pp\Big(((L'_{(1,2)}(t),\dots,L'_{(n-1,n)}(t)),\;t\in[0,T])\in A\Big)\geq\frac{1}{2},\\
&&\pp\Big(((X'_{(2)}(t)-X'_{(1)}(t),\dots,X'_{(n)}(t)-X'_{(n-1)}(t)),\;t\in[0,T])\in B\Big)\geq\frac{1}{2}
\end{eqnarray}
and for any $r>0$ set  
\begin{eqnarray}
A_r=\{h\in C([0,T],\rr^{n-1})|\;\inf_{\widetilde{h}\in A}\|\widetilde{h}-h\|_{T,2}\leq r\}, \\
B_r=\{h\in C([0,T],\rr^{n-1})|\;\inf_{\widetilde{h}\in B}\|\widetilde{h}-h\|_{T,2}\leq r\}. 
\end{eqnarray}
Then for all $r_1\geq 2^4\sqrt{\frac{2(K-1)^5 T \log 2}{n-1}}$ and $r_2\geq 2^4 3\sqrt{\frac{2(K-1)^5 T \log 2}{n-1}}$ it holds
\begin{eqnarray*}
&&\pp\Big(((L'_{(1,2)}(t),\dots,L'_{(n-1,n)}(t)),\;t\in[0,T])\notin A_{r_1}\Big)\leq \exp\Big(-\frac{r_1^2(n-1)}{2^9 (K-1)^5 T}\Big),\\
&&\pp\Big(((X'_{(2)}(t)-X'_{(1)}(t),\dots,X'_{(n)}(t)-X'_{(n-1)}(t)),\;t\in[0,T])\notin B_{r_2}\Big)
\leq \exp\Big(-\frac{r_2^2(n-1)}{2^9 3^2(K-1)^5 T}\Big). 
\end{eqnarray*}}
{\it Proof.} From the considerations in step 1 of the proof of Theorem \ref{slopetheorem} we see that the process 
\begin{eqnarray} \label{lproc}
(L'_{(1,2)}(t),\dots,L'_{(n-1,n)}(t)),\;t\in[0,T]
\end{eqnarray}
can be obtained by applying the map $\Phi^{Q^{(K-1)}}_L$ and then the canonical projection of $C([0,T],\rr^{K-1})$ onto $C([0,T],\rr^{n-1})$ to the process
\begin{eqnarray}\label{reflproc}
\Big((\delta_2-\delta_1)t+\beta'_2(t)-\beta'_1(t),\dots,(\delta_K-\delta_{K-1})t+\beta'_K(t)-\beta'_{K-1}(t)\Big),\;t\in[0,T].
\end{eqnarray}
Hereby, the Brownian motions $\beta'_1,\dots,\beta'_K$ are defined analogously to the Brownian motions $\beta_1,\dots,\beta_K$ in the proof of Theorem \ref{slopetheorem}. There, we have seen that the process in \eqref{reflproc} satisfies a QTCI with respect to the norm $\|.\|_{T,2}$ with the constant $16(K-1)^{-1}T$. Combining Lemma \ref{lemmapushforward} and Lemma \ref{lemma_lipschitz}, we conclude that a QTCI with the constant 
\begin{eqnarray}
\Big(\sqrt{\frac{K-1}{n-1}}\Big)^2\cdot\Big(2(K-1)^{\frac{5}{2}}\Big)^2\cdot\frac{16 T}{K-1}=\frac{2^6(K-1)^5T}{n-1} 
\end{eqnarray}
applies to the process in (\ref{lproc}) with respect to the norm $\|.\|_{T,2}$. Hence, from Theorem \ref{preq} we obtain the first of the two claimed concentration of measure results. 

\bigskip

To show the second one, we recall from step 1 of the proof of Theorem \ref{slopetheorem} that the process
\begin{eqnarray}\label{yproc}
(X'_{(2)}(t)-X'_{(1)}(t),\dots,X'_{(n)}(t)-X'_{(n-1)}(t)),\;t\in[0,T]
\end{eqnarray}
is the image of the process in (\ref{reflproc}) under the successive application of the map $\Phi^{Q^{(K-1)}}_R$ and the canonical projection of $C([0,T],\rr^{K-1})$ onto $C([0,T],\rr^{n-1})$.

Moreover, we can rewrite \eqref{skorokhod} as
\begin{eqnarray}
(\Phi^{Q^{(K-1)}}_R(h))(t)=h(t)+(I^{(K-1)}-Q^{(K-1)})((\Phi^{Q^{(K-1)}}_L(h))(t)),\quad t\in[0,T]
\end{eqnarray}
for all $h\in C([0,T],\rr^{K-1})$, where $I^{(K-1)}$ is the $(K-1)\times(K-1)$ identity matrix. It follows that for all $h_1,h_2\in C([0,T],\rr^{K-1})$ one has the estimates
\begin{eqnarray*}
&&\Big\|(\Phi^{Q^{(K-1)}}_R(h_2))-(\Phi^{Q^{(K-1)}}_R(h_1))\Big\|_{T,2}\\
&&\leq \|h_2-h_1\|_{T,2}+\Big\|(I^{(K-1)}-Q^{(K-1)})((\Phi^{Q^{(K-1)}}_L(h_2))(.)-(\Phi^{Q^{(K-1)}}_L(h_1))(.))\Big\|_{T,2}\\
&&\leq \|h_2-h_1\|_{T,2}+\frac{3\sqrt{2}}{2}\cdot\Big\|(\Phi^{Q^{(K-1)}}_L(h_2))-(\Phi^{Q^{(K-1)}}_L(h_1))\Big\|_{T,2}\\
&&\leq \Big(1+3\sqrt{2}(K-1)^{\frac{5}{2}}\Big)\|h_2-h_1\|_{T,2}. 
\end{eqnarray*} 
In the second inequality we have combined the fact that the matrix $I^{(K-1)}-Q^{(K-1)}$ is tridiagonal with the elementary inequality $(a_1+a_2+a_3)^2\leq 3(a_1^2+a_2^2+a_3^2)$, $a_1,a_2,a_3\in\rr$. The third inequality is a consequence of Lemma \ref{lemma_lipschitz}.
 
Now, it follows from Lemma \ref{lemmapushforward} that the process in (\ref{yproc}) satisfies a QTCI with respect to the norm $\|.\|_{T,2}$ with the constant 
\begin{eqnarray}
\left(\sqrt{\frac{K-1}{n-1}}\right)^2\cdot\Big(1+3\sqrt{2}(K-1)^{\frac{5}{2}}\Big)^2\cdot\frac{16 T}{K-1}\leq\frac{2^6 3^2(K-1)^5 T}{n-1}.
\end{eqnarray}
The second claim of the corollary is a consequence of this and Theorem \ref{preq}.\ep 

\bigskip

The last ingredient in the proof of Theorem \ref{thm_infinite} is an estimate on how fast the stopping times $0=\tau_N\leq\tau_{N+1}\leq\dots$ in the proof of Lemma \ref{lemmastopping} grow to infinity in terms of the initial positions of the particles. 
{\lemma\label{taulemma} Let the Assumption 1.1 be satisfied with a constant $c>0$. Then for all natural numbers $m\geq\frac{\max_{j=1,\dots,M-1}|\delta_j-\delta_M| T}{c}+1$ ($=\frac{\Delta T}{c}+1$) one has the inequality
\begin{eqnarray}\label{boundstopping}
\pp(\tau_{N+m}\leq T)\leq \frac{NT}{c(cm-c-\Delta T)}\cdot\exp\Big(-\frac{1}{2T}(cm-c-\Delta T)^2\Big).
\end{eqnarray}
In particular, there exists a constant $C(c,M,n,T,\Delta)>0$ independent of $m$ such that
\begin{eqnarray}\label{cbound}
\pp(\tau_{N+m}\leq T)\leq C(c,M,n,T,\Delta)e^{-\frac{c^2}{3T} m^2},\qquad m\in\nn. 
\end{eqnarray}}
{\it Proof.} We fix a natural number $m$ as in the first statement of the lemma and note that on the event $\{\tau_{N+m}\leq T\}$ there exist numbers $1\leq i\leq N$ and $j\geq N+m$ such that the particle, which was the $i$-th from the left in the initial particle configuration, appears on the right or at the same position as the particle, which was the $j$-th from the left in the initial particle configuration, at a time $t\in[0,T]$. Using this observation, the union bound and the definition of $\Delta$ (see the statement of Theorem \ref{thm_infinite}), one has the chain of inequalites  
\begin{eqnarray*}
\pp(\tau_{N+m}\leq T)&\leq&\sum_{i=1}^N\sum_{j=N+m}^\infty \pp(\sup_{0\leq t\leq T} (W_i(t)-W_j(t))\geq -\Delta T+X_j(0)-X_i(0))\\
&\leq& N \sum_{j=N+m}^\infty\pp(\sup_{0\leq t\leq T} (W_1(t)-W_j(t))\geq -\Delta T+X_j(0)-X_N(0)).
\end{eqnarray*}

From Bernstein's inequality for Brownian motion (see page 145 in \cite{ry}), Assumption \ref{assumption} and the assumption $m\geq\frac{\Delta c}{T}+1$ it follows that the latter expression can be bounded further by
\begin{eqnarray*}
&&N\sum_{j=N+m}^\infty \exp\Big(-\frac{(X_j(0)-X_N(0)-\Delta T)^2}{2T}\Big)
\leq N\sum_{k=m}^\infty \exp\Big(-\frac{(ck-\Delta T)^2}{2T}\Big)\\
&&=N\sum_{k=m}^\infty \exp\Big(-\frac{c^2}{2T}\Big(k-\frac{\Delta T}{c}\Big)^2\Big)
\leq N\int_{m-1}^\infty \exp\Big(-\frac{c^2}{2T}\Big(y-\frac{\Delta T}{c}\Big)^2\Big)\;dy.
\end{eqnarray*}

Next, we note that the latter integral is equal to the probability that a standard normal random variable exceeds 
$\frac{m-1-\frac{\Delta T}{c}}{\frac{\sqrt{T}}{c}}$ multiplied by $\sqrt{2\pi\frac{T}{c^2}}$. Using this and the standard estimate 
\begin{eqnarray}
\int_y^\infty e^{-\frac{z^2}{2}}\;dz\leq\frac{1}{y}e^{-\frac{y^2}{2}},\quad y>0
\end{eqnarray}
one ends up with the first statement of the lemma. 

Finally, to see (\ref{cbound}), it suffices to observe that the argument of the exponential function on the right-hand side of the inequality \eqref{boundstopping} is a quadratic polynomial in $m$, in which the coefficient of $m^2$ is given by $\frac{c^2}{2T}$. \ep 

\bigskip

We can now prove the following refined version of Theorem \ref{thm_infinite}.

{\prop Let the sets $A$, $B$, $A_r$, $r>0$ and $B_r$, $r>0$ be defined as in Theorem \ref{thm_infinite} and let the constant $C(c,M,n,T,\Delta)$ be as in Lemma \ref{taulemma}. Moreover, let $m_1,m_2\in\nn$ be such that for all natural numbers $m\geq m_1$ (or $m\geq m_2$) the value of $C(c,M,n,T,\Delta)e^{-\frac{c^2}{3T}m^2}$ is less or equal to the difference between the left-hand and the right-hand side of the inequality (\ref{Acond}) (or (\ref{Bcond}), respectively). Also, define $C_1$ and $C_2$ to be the smallest positive real numbers such that 
\begin{eqnarray}
&&\left(\frac{3(n-1)}{2^9 c^2}\right)^{\frac{1}{7}} C_1^{\frac{2}{7}}\geq m_1,\label{C1cond1}\\
&&\forall r\geq C_1:\;\;r\geq 2^4\sqrt{\frac{2(N+m_1(r)-1)^5 T \log 2}{n-1}},\label{C1cond2}\\
&&\left(\frac{n-1}{2^9 3 c^2}\right)^{\frac{1}{7}} C_2^{\frac{2}{7}}\geq m_2,\label{C2cond1}\\
&&r\geq C_2:\;\;r\geq 2^4\cdot3\sqrt{\frac{2(N+m_1(r)-1)^5 T \log 2}{n-1}},\label{C2cond2}
\end{eqnarray}
where we have set $m_1(r)=\left(\frac{3(n-1)}{2^9 c^2}\right)^{\frac{1}{7}} r^{\frac{2}{7}}$ and 
$m_2(r)=\left(\frac{n-1}{2^9 3 c^2}\right)^{\frac{1}{7}} r^{\frac{2}{7}}$. Then there exist positive constants $C_3$, $C_4$ depending on $c$, $\Delta$, $M$, $n$, $T$ and the value on the left-hand side of (\ref{Acond}) and (\ref{Bcond}), respectively, such that for all $r_1\geq C_1$ and $r_2\geq C_2$ it holds
\begin{eqnarray*}
\pp\Big(((L_{(1,2)}(t),\dots,L_{(n-1,n)}(t)),t\in[0,T])\notin A_{r_1}\Big)
\leq C_3\exp\Big(-r_1^{\frac{4}{7}}\cdot\frac{(n-1)^{\frac{2}{7}}c^{\frac{10}{7}}}{2^{\frac{18}{7}}3^{\frac{5}{7}}T}\Big), \label{Arineq2}\\
\pp\Big(((X_{(2)}(t)-X_{(1)}(t),\dots,X_{(n)}(t)-X_{(n-1)}(t)),t\in[0,T])\notin B_{r_2}\Big)
\leq C_4\exp\Big(-r_2^{\frac{4}{7}}\cdot\frac{(n-1)^{\frac{2}{7}}c^{\frac{10}{7}}}{2^{\frac{18}{7}} 3^{\frac{9}{7}}T}\Big). \label{Brineq2}
\end{eqnarray*}}
{\it Proof.} Since the way of proof is the same for both inequalities, we only provide the proof of the first one. To this end, we fix an $r_1$ as in the statement of the proposition and set 
\begin{eqnarray}
&&\widetilde{m}_1=\Big(\frac{3(n-1)}{2^9 c^2}\Big)^{\frac{1}{7}} r_1^{\frac{2}{7}},\\
&&K_1=N+\widetilde{m}_1. 
\end{eqnarray}
We assume from now on that $r_1$ is such that $\tilde{m}_1$ is an integer. If this is not the case, one merely needs to replace $\tilde{m}_1$ by the smallest integer which is larger than $\tilde{m}_1$. Moreover, we note that the inequalities $r_1\geq C_1$ and \eqref{C1cond1} imply $\widetilde{m}_1\geq m_1$. 

Next, we recall the definition of the stopping time $\tau_{K_1}$ and observe
\begin{eqnarray*}
&&\pp\Big(((L_{(1,2)}(t),\dots,L_{(n-1,n)}(t)),t\in[0,T])\notin A_{r_1}\Big)\\
&&\leq\pp\Big(((L_{(1,2)}(t),\dots,L_{(n-1,n)}(t)),t\in[0,T])\notin A_{r_1},\tau_{K_1}\geq T\Big)+\pp(\tau_{N+\widetilde{m}_1}\leq T)\\
&&\leq\pp\Big(((L_{(1,2)}(t),\dots,L_{(n-1,n)}(t)),t\in[0,T])\notin A_{r_1},\tau_{K_1}\geq T\Big)
+C(c,M,n,T,\Delta)e^{-\frac{c^2}{3T}\widetilde{m}_1^2},
\end{eqnarray*}
where the last inequality is a consequence of Lemma \ref{taulemma}. 

To obtain an upper bound on the first summand, which we call term (*), we remark first that for all triples of indices $1\leq i<j<k$ it holds
\begin{eqnarray}\label{triplecoll}
\pp(\exists t\in[0,T]:\;X_{(i)}(t)=X_{(j)}(t)=X_{(k)}(t))=0.
\end{eqnarray}
Indeed, arguing as in the proof of Lemma \ref{lemmastopping}, but replacing $N$ by $\max(k,M)$, we deduce the existence of stopping times $0=\widetilde{\tau}_{\max(k,M)}\leq \widetilde{\tau}_{\max(k,M)+1}\leq\dots$ tending to infinity almost surely and such that for each $l\in\nn$ the dynamics of the $k$ leftmost particles, stopped at $\widetilde{\tau}_{\max(k,M)+l}$, is given by the dynamics of the $k$ leftmost particles in a finite particle system as in (\ref{ranksde}), stopped at $\widetilde{\tau}_{\max(k,M)+l}$. Hence, (\ref{triplecoll}) is a consequence of the considerations in section 2.2 of \cite{ik}, Proposition 1 in \cite{ik} and the fact that a countable union of $\pp$-null sets is $\pp$-null set. 

Next, we recall from the proof of Lemma \ref{lemmastopping} that for each $m\in\nn$ the paths of the process $X_{(1)}(t\wedge\tau_{N+m}),\dots,X_{(n)}(t\wedge\tau_{N+m})$, $t\in[0,T]$ are the paths of the $n$ leftmost particles in a particle system as in \eqref{ranksde} with $|\Lambda_{N+m}|$ particles, stopped at $\tau_{N+m}$. Moreover, (\ref{triplecoll}) shows that $|\Lambda_{N+m}|=N+m$. Thus, on the event $\{\tau_{K_1}\geq T\}$ the paths of the process 
\begin{eqnarray*}
(L_{(1,2)}(t),\dots,L_{(n-1,n)}(t)),t\in[0,T]
\end{eqnarray*}    
can be written as the composition of the map $\Phi^{Q^{(K_1-1)}}_L$ with the canonical projection of $C([0,T],\rr^{K_1-1})$ onto $C([0,T],\rr^{n-1})$ applied to the paths of the process
\begin{eqnarray*}
\Big((\delta_2-\delta_1)t+\beta^{(\widetilde{m}_1)}_2(t)-\beta^{(\widetilde{m}_1)}_1(t),\dots,
(\delta_{K_1}-\delta_{K_1-1})t+\beta^{(\widetilde{m}_1)}_{K_1}(t)-\beta^{(\widetilde{m}_1)}_{K_1-1}(t)\Big),\;t\in[0,T],
\end{eqnarray*}
where $\beta^{(\widetilde{m}_1)}_1,\dots,\beta^{(\widetilde{m}_1)}_{K_1}$ are i.i.d. standard Brownian motions defined in Lemma \ref{lemmastopping}. Hence, following the proof of the first statement of Corollary \ref{corfinite} and using \eqref{C1cond2}, we conclude that term (*) is bounded from above by $\exp\Big(-\frac{r_1^2(n-1)}{2^9 (K_1-1)^5 T}\Big)$. 

All in all, we have shown that $\pp\Big(((L_{(1,2)}(t),\dots,L_{(n-1,n)}(t)),t\in[0,T])\notin A_{r_1}\Big)$ is bounded from above by
\begin{eqnarray*}
\exp\Big(-\frac{r_1^2(n-1)}{2^9 (K_1-1)^5 T}\Big)+C(c,M,n,T,\Delta)e^{-\frac{c^2}{3T}\widetilde{m}_1^2}.
\end{eqnarray*}
Plugging in the values of $\widetilde{m}_1$ and $K_1$ one observes that the leading order term in the variable $r_1$ is the same for both exponents on the right-hand side of the latter inequality and equals to $-r_1^{\frac{4}{7}}\cdot\frac{(n-1)^{\frac{2}{7}}c^{\frac{10}{7}}}{2^{\frac{18}{7}}3^{\frac{5}{7}}T}$. This immediately finishes the proof of the proposition and completes the proof of Theorem \ref{thm_infinite}. \ep


\begin{thebibliography}{20}

\bibitem{AS}
\textsc{Abramowitz, M.} and \textsc{Stegun, I.} (eds.) (1984)
\newblock {\em Pocketbook of mathematical functions.} 
Abridged edition of Handbook of mathematical functions. 
Material selected by \textsc{Michael Danos} and \textsc{Johann Rafelski}. Verlag Harri Deutsch, Thun, 1984. 

\bibitem{AA}
\textsc{Arguin, L.~-P.} and \textsc{Aizenman, M.} (2009). 
On the structure of quasi-stationary competing particles systems. 
\textit{Ann. Probab.} \textbf{37} 1080-1113.

\bibitem{arratia83}
\textsc{Arratia, R.} (1983).
The motion of a tagged particle in the simple symmetric exclusion system on ${\bf {Z}}$.
\textit{Ann. Probab.} \textbf{11} 362-373.

\bibitem{atlasmodel}
\textsc{Banner, A.}, \textsc{Fernholz, R.} and \textsc{Karatzas, I.} (2005).
Atlas models of equity markets.
\textit{Ann. Appl. Probab.} \textbf{15} 2296-2330.

\bibitem{BG}
\textsc{Banner, A.} and \textsc{Ghomrasni, R.} (2008).
Local times of ranked continuous semimartingales. 
\textit{Stochastic Process. Appl.} \textbf{118} 1244-1253.

\bibitem{BGL}
\textsc{Bobkov, S.}, \textsc{Gentil, I.} and \textsc{Ledoux, M.} (2001).
Hypercontractivity of Hamilton-Jacobi equations. \textit{J. Math. Pure Appl.} \textbf{80} 669-696.

\bibitem{chatpal}
\textsc{Chatterjee, S.} and \textsc{Pal, S.} (2010).
\newblock A phase transition behavior for Brownian motions interacting through their ranks. 
\textit{Probab. Theory Related Fields} \textbf{147} 123-159. 

\bibitem{chatpal2}
\textsc{Chatterjee, S.} and \textsc{Pal, S.} (2008).
\newblock A combinatorial analysis of interacting diffusions. 
To appear in \textit{J. Theor. Probab.}.

\bibitem{demasiferrari}
\textsc{De Masi, A.} and \textsc{Ferrari, P.} (2002).
Flux fluctuations in the one dimensional nearest neighbors symmetric simple exclusion process.
\textit{J. Stat. Phys.} \textbf{107} 677-683.  

\bibitem{D}
\textsc{Dembo, A.} (1997)
Information inequalities and concentration of measures. 
\textit{Ann. Probab.} \textbf{25} 927-939.

\bibitem{DZ}
\textsc{Dembo, A.} and \textsc{Zeitouni, O.} (1996)
Transportation approach to some concentration inequalities in product spaces. 
\textit{Electron. Comm. Probab.} \textbf{1}, 83-90.

\bibitem{DGW}
\textsc{Djellout, H.}, \textsc{Guillin, A.} and \textsc{Wu, L.} (2004).
Transportation cost-information inequalities and applications to random dynamical systems and diffusions. 
\textit{Ann. Probab.} \textbf{32} 2702-2732.

\bibitem{DI}
\textsc{Dupuis, P.} and \textsc{Ishii, H.} (1991).
On Lipschitz continuity of the solution mapping to the Skorokhod problem, with applications. 
\textit{Stochastics} \textbf{35}, 31--62.

\bibitem{DR1}
\textsc{Dupuis, P.} and \textsc{Ramanan, K.} (1999).
Convex Duality and the Skorokhod Problem. I. 
\textit{Probab. Theory Related Fields} \textbf{115} 153-195.

\bibitem{DR2}
\textsc{Dupuis, P.} and \textsc{Ramanan, K.} (1999).
Convex Duality and the Skorokhod Problem. II. 
\textit{Probab. Theory Related Fields} \textbf{115} 197-236.

\bibitem{FS05}
\textsc{Fang, S.} and \textsc{Shao, J.} (2005). 
Transportation cost inequalities on path and loop groups.
\textit{J. Funct. Anal.} \textbf{218} 293-317.

\bibitem{FS07}
\textsc{Fang, S.} and \textsc{Shao, J.} (2007). 
Optimal transport maps for Monge-Kantorovich problem on loop
groups. 
\textit{J. Funct. Anal.} \textbf{248} 225-257.

\bibitem{FWW}
\textsc{Fang, S.}, \textsc{Wang, F.~-Y.}, and \textsc{Wu, B.} (2008). 
Transportation-cost inequality on path spaces with uniform  distance.  
\textit{Stochastic Process. Appl.}  \textbf{118} 2181-2197.

\bibitem{fern02}
	\textsc{Fernholz, R.} (2002).
	\newblock {\em Stochastic Portfolio Theory}.
	\newblock Springer, New York.

\bibitem{FK}
\textsc{Fernholz, R.} and \textsc{Karatzas, I.} (2009). 
Stochastic Portfolio Theory: A survey. In \textit{Handbook of Numerical Analysis: Mathematical Modeling and Numerical Methods in Finance}. Elsevier Publishing Company BV, Amsterdam, 89--168.

\bibitem{FU}
\textsc{Feyel, D.} and \textsc{\"Ust\"unel, A.~S.} (2004). 
The Monge-Kantorovitch problem and Monge-Amp\`ere equation on Wiener space. 
To appear in \textit{Probab. Theory Related Fields}.

\bibitem{GW}
\textsc{Gourcy, M.} and \textsc{Wu, L.} (2006).
Logarithmic Sobolev inequalities of diffusions for the $L^2$ metric. 
\textit{Potential Anal.}  \textbf{25} 77-102.

\bibitem{harris65}
\textsc{Harris, T. E.} (1965).
\newblock Diffusion with ``collisions'' between particles.
\newblock {\em J. Appl. Probab.} \textbf{2} 323-338.

\bibitem{hr} \textsc{Harrison, J. M.} and \textsc{Reiman, M. I.} (1981). 
Reflected Brownian motion in an orthant. 
\textit{Ann. Probab.} \textbf{9} 302-308.

\bibitem{HP}
\textsc{Houdr\'e, C.} and \textsc{Privault, N.} (2002).
Concentration and deviation inequalities in infinite dimensions via covariance representations.
\textit{Bernoulli} \textbf{8} 697-720.

\bibitem{ik} \textsc{Ichiba, T.} and \textsc{Karatzas, I.} (2009). 
On collisions of Brownian particles. 
Preprint available at http://arxiv.org/abs/0810.2149v2. To appear in \textit{Ann. Appl. Probab.} 

\bibitem{IPBKF} 
\textsc{Ichiba, T.}, \textsc{Papathanakos, V.}, \textsc{Banner, A.}, \textsc{Karatzas, I.}, and \textsc{Fernholz, R.} (2010). 
Hybrid Atlas Models. 
Preprint available at http://arxiv.org/abs/0909.0065. To appear in \textit{Ann. Appl. Probab.}

\bibitem{joumal}
\textsc{Jourdain, B.} and \textsc{Malrieu, F.} (2008).
\newblock Propagation of chaos and Poincar\'{e} inequalities for a system of particles interacting through their cdf.
\newblock \textit{Ann. Appl. Probab.} \textbf{18} 1706-1736. 

\bibitem{L}
\textsc{Ledoux, M.} (2001)
\textit{The Concentration of Measure Phenomenon.}
Mathematical Surveys and Monographs \textbf{89}. American Mathematical Society.

\bibitem{M1} 
\textsc{Marton, K.} (1996).
Bounding $\bar d$-distance by information divergence: a method to prove measure concentration. 
\textit{Ann. Probab.} \textbf{24} 857-866.

\bibitem{M2}
\textsc{Marton, K.} (1997).
A measure concentration inequality for contracting Markov chains.
\textit{Geom. Funct. Anal.} \textbf{6} 556-571. 

\bibitem{M3}
\textsc{Marton, K.} (1998).
Mesure concentration for a class of random processes.
\textit{Probab. Theory Related Fields} \textbf{110} 427-439.

\bibitem{sheppmckean}
\textsc{McKean, H. P.} and \textsc{Shepp, L.} (2005).
\newblock The advantage of capitalism vs. Socialism depends on the criterion.
\newblock Available at {\em www.emis.de/journals/ZPOMI/v328/p160.ps.gz}.

\bibitem{NV}
\textsc{Nourdin, I.} and \textsc{Viens, F.~G.} (2009).
Density formula and concentration inequalities with Malliavin calculus.
\textit{Electron. J. Probab.} \textbf{14} 2287-2309.

\bibitem{OV}
\textsc{Otto, F.} and \textsc{Villani, C.} (2000).
Generalization of an inequality by Talagrand and links with the logarithmic Sobolev inequality. 
\textit{J. Funct. Anal.} \textbf{173} 361-400.

\bibitem{pa} \textsc{Pal, S.} (2010). Concentration for multidimensional diffusions and their boundary local times. Accepted (modulo minor revision) in \textit{Probab. Theory Related Fields}. Preprint available at http://arxiv.org/abs/1005.2217. 

\bibitem{pp} \textsc{Pal, S.} and \textsc{Pitman J.} (2008). One-dimensional Brownian particle systems with rank-dependent drifts. 
\textit{Ann. Appl. Probab.} \textbf{18} 2179-2207.

\bibitem{ry} \textsc{Revuz, D.} and \textsc{Yor, M.} {\em Continuous martingales and Brownian motion}. (1999). Springer-Verlag, Berlin. 3rd ed.  

\bibitem{rostvares}
\textsc{Rost, H.} and \textsc{Vares, M. E.} (1985).
\newblock Hydrodynamics of a one-dimensional nearest neighbor model.
\newblock{\em Contemp. Math.} \textbf{41} 329-342.

\bibitem{ruzaizenman}
\textsc{Ruzmaikina, A.} and \textsc{Aizenman, M.} (2005).
Characterization of invariant measures at the leading edge for competing particle systems.
\textit{Ann. Probab.}, \textbf{33} (1), 82-113.

\bibitem{shkol}
\textsc{Shkolnikov, M.} (2009). Competing Particle Systems Evolving by I.I.D. Increments. 
\textit{Electron. J. Probab.} \textbf{14} 728-751. 

\bibitem{sh} \textsc{Shkolnikov, M.} (2010). Competing particle systems evolving by interacting Levy processes. Preprint available at http://arxiv.org/abs/1002.2811. To appear in \textit{Ann. Appl. Probab.}

\bibitem{sh2} \textsc{Shkolnikov, M.} (2010). Large systems of diffusions interacting through their ranks. Preprint available at http://arxiv.org/abs/1008.4611. 

\bibitem{Sw} \textsc{Swanson, J.} (2007). Weak convergence of the scaled median of independent Brownian motions. \textit{Probab. Theory Related Fields} \textbf{138} 269-304.

\bibitem{sznitman86}
\textsc{Sznitman, A. S.} (1986).
\newblock A propagation of chaos result for {B}urgers' equation.
\newblock {\em Probab. Theory Related Fields} \textbf{71} 581-613.

\bibitem{sznitman91}
\textsc{Sznitman, A. S.} (1991).
\newblock Topics in propagation of chaos.
\newblock In {\em \'Ecole d'\'Et\'e de Probabilit\'es de Saint-Flour
  XIX---1989}, 165--251. Springer, Berlin.

\bibitem{T96}
\textsc{Talagrand, M.} (1996). 
A new look at independence. 
\textit{Ann. Probab.} \textbf{24} 1-34. 

\bibitem{T96b}
\textsc{Talagrand, M.} (1996). 
Transportation cost for Gaussian and other product measures. 
\textit{Geom. Funct. Anal.} \textbf{6} 587-600.

\bibitem{U}
\textsc{\"Ust\"unel, A. S.} (2010). 
Transportation cost inequalities for diffusions under uniform distance. 
Preprint available at http://arxiv.org/abs/1009.5251.

\bibitem{W02}
\textsc{Wang, F.~-Y.} (2002). 
Transportation cost inequalities on path spaces over Riemannian manifolds.
\textit{Illinois J. Math.} \textbf{46} 1197-1206.

\bibitem{W08}
\textsc{Wang, F.~-Y.} (2008). 
Generalized transportation-cost inequalities and applications. 
\textit{Potential Anal.} \textbf{28} 321-334.

\bibitem{WZ}
\textsc{Wu, L.} and \textsc{Zhang, Z.} (2004). 
Talagrand's $T_2$-transportation inequality w.r.t. a uniform metric for diffusions.  
\textit{Acta Math. Appl. Sin. Engl. Ser.}  \textbf{20} 357-364.

\end{thebibliography}
\end{document}